\documentclass[12pt,reqno]{amsart}
\usepackage{etex}
\usepackage{amssymb}
\usepackage{mathrsfs}
\usepackage{palatino}
\usepackage{bbm}
\usepackage{pdflscape}
\usepackage{mathrsfs}
\usepackage{cases}
\usepackage{epic}
\usepackage{amsfonts}
\usepackage{graphicx}
\usepackage{amsmath}
\usepackage{amssymb, upgreek}
\usepackage{bm}
\usepackage{latexsym,todonotes}
\usepackage{pdflscape}
\usepackage[all]{xypic}
\usepackage{color}
\usepackage{colordvi}
\usepackage{multicol}
\usepackage[normalem]{ulem}
\usepackage[linktocpage=true]{hyperref}
\usepackage{hyperref}
\hypersetup{colorlinks,linkcolor=blue,urlcolor=cyan,citecolor=red}
\usepackage{pictex}
\usepackage{amscd}
\usepackage{latexsym}
\usepackage{amssymb}
\usepackage{eucal}
\usepackage{eufrak}
\usepackage{verbatim}

\newcommand{\NN}{\mathbb{N}}

\newcommand{\fg}{\mathfrak{g}}

\newcommand{\kk}{\mathbbm{k}}

\newcommand{\gl}{\mathfrak{gl}}

\newcommand{\fS}{\mathfrak{S}}

\newcommand{\fo}{\mathfrak{o}}

\newcommand{\xon}{{\rm X}(\mathfrak{o}_{N})}
\newcommand{\yon}{{\rm Y}(\mathfrak{o}_{N})}

\newcommand{\tka}{\tilde{\kappa}}
\newcommand{\Qpm}{Q^{\pm}}

\DeclareMathOperator{\ad}{ad}

\DeclareMathOperator{\Char}{char}

\DeclareMathOperator{\gr}{gr}

\DeclareMathOperator{\X}{X}
\DeclareMathOperator{\Y}{Y}

\DeclareMathOperator{\HC}{HC}

\DeclareMathOperator{\U}{U}

\numberwithin{equation}{section}
\newtheorem{Theorem}{Theorem}[section]
\newtheorem{Lemma}[Theorem]{Lemma}
\newtheorem{Corollary}[Theorem]{Corollary}
\newtheorem{Proposition}[Theorem]{Proposition}
\newtheorem{Remark}[Theorem]{Remark}

\theoremstyle{Theorem}

\newtheorem*{thm*}{Theorem}
\newtheorem*{thm**}{Corollary}
\newtheorem*{thm***}{Theorem B}
\theoremstyle{remark}

\numberwithin{equation}{section}

\begin{document}
\title[Orthogonal]{Modular orthogonal Yangians}
\author[Hao Chang \lowercase{and} Hongmei Hu]{Hao Chang \lowercase{and} Hongmei Hu*}
\address[Hao Chang]{School of Mathematics and Statistics and Hubei Key Laboratory of Mathematical Sciences, Central China Normal University, Wuhan 430079, China}
\email{chang@ccnu.edu.cn}
\address[Hongmei Hu]{Department of Mathematics, Shanghai Maritime University, Shanghai 201306, China}
\email{hmhu@shmtu.edu.cn}
\date{\today}
\thanks{* Corresponding author.}


\begin{abstract}
We study the (extended) orthogonal Yangians associated to the Lie algebras types $B$ and $D$ over a field of positive characteristic.
We define the $p$-center for the Yangians and obtain an explicit description of the center in terms of Drinfeld generators,
showing that the center is generated by its Harish-Chandra center together with a large $p$-center.
\end{abstract}
\maketitle

\section{Introduction}
For each simple Lie algebra $\fg$ over the field of complex numbers,
the corresponding Yangian $\Y(\fg)$ was defined by Drinfeld in \cite{D85}
as a cannonical deformation of the universal enveloping algerba $\U(\fg[x])$
for the {\it current Lie algebra} $\fg[x]$.
In \cite{D88}, Drinfeld gave a new presentation which is now often referred to as the {\it Drinfeld presentation}.
The Yangian $\Y_n$ associated to the general Lie algebra $\gl_n$ was earlier considered in the work
of mathematical physicists from St. Petersburg \cite{TF79}.
It is an associative algebra whose defining relations can be written in a specific matrix form,
which is called the {\em RTT relation}; see e.g. \cite{MNO96}.
In \cite{AACFR03} (see also \cite{AMR06}),
the RTT-presentation of the Yangian associated with $B,C$ or $D$ Lie algebras was studied.
Regarding for related topics and further applications of Yangians,
the reader is referred to \cite{Mol07} and the references cited therein.

An explicit isomorphism between the RTT and Drinfeld presentations of the Yangian $\Y_n$
is constructed with the use of the {\it Gauss decomposition} (or {\it quasideterminants}) of the generator matrix.
Complete proofs were given by Brundan and Kleshchev \cite{BK05},
see also \cite[Chapter 3]{Mol07}.
The same approach for Yangians of type BCD is more challenging and was accomplished by Jing, Liu and Molev \cite{JLM18} ,
while a different method to establish isomorphisms was developed in \cite{GRW19}.

In \cite{BT18}, Brundan and Topley developed the theory of the Yangian $\Y_n$ over a field of positive characteristic.
One of the key features which differs from characteristic zero is the existence of a large central subalgebra $Z_p(\Y_n)$,
called the {\it $p$-center}.
Very recently, the first author jointly with Hu and Topley \cite{CHT25} studied the modular representations of the Yangian $\Y_2$.
Using the $p$-center, they classified the finite-dimensional irreducible modules for the {\it restricted Yangian} $\Y_2^{[p]}$.
This leads to a description of the certain non-restricted representations of the higher rank general linear Lie algebras.
From this one needs precise information about the structure of the center of the Yangian in positive characteristic.

Our goal in this paper is to give a description of the center $Z(\yon)$ of the modular orthogonal Yangian $\yon$.
Let us shortly explain our approach.
Actually, we will mainly work with the {\it extended orthogonal Yangian} $\xon$.
We first define it over an algebraically closed field $\kk$ of positive characteristic to be the associative algebra by the usual RTT-relation (\eqref{def relations RTT}).
Our main methods are similar strategy to \cite{BT18}, see also \cite{CH23} for the type $A$ super Yangians.
To describe the center $Z(\xon)$,
it was easier to work with the Drinfeld presentation rather than the RTT-presentation of $\xon$.
Under the assumption the characteristic $\Char\kk=:p>3$,
it was observed by Brundan-Topley (\cite[Theorem 4.3]{BT18}) that the Drinfeld presentation of the modular Yangian $\Y_n$
is exactly the same as the one obtained by Brundan-Kleshchev over the complex field (see \cite[Theorem 5.2]{BK05}).
It remains to be true in our setting, we obtain the Drinfeld presentation of the modular orthogonal Yangians
by using the results of Jing-Liu-Molev \cite{JLM18} (see Theorem \ref{Theorem: Drinfeld presentation}).
However, the proof of the Serre relations there
does not work in our situation.
We establish the Serre relations by adapting the methods of \cite{Lev93,GNW18} to positive characteristic (see Appendix \ref{section name:appendix}).
Moreover, as the Drinfeld presentation involves more relations,
the calculations become more complicated than that in \cite{BT18, CH23}.
We need to put extra effort to obtain more formulae (Lemma \ref{useful lemma 1}) for getting the central elements.

Next, we take advantage of the {\it transposition} and {\it permutation} automorphisms of $\xon$,
These automorphisms reduces the verification of all the relations to checking the simple generating elements,
so that we may derive various relations by direct computations.
Notice that the root system contains roots of two distinct lengths for type $B$ Lie algebras,
we must deal with them separately (see Section \ref{subsec: anto}).
Finally we show that we have found enough generating elements for the center $\xon$.

In characteristic zero, the center $Z(\xon)$ was determined in \cite{AACFR03, AMR06}.
It still makes sense when $\Char\kk =p>0$.
We call it the {\it Harish-Chandra center} of $\xon$ and denote it by $Z_{\HC}(\xon)$.
In positive characteristic, the orthogonal current Lie algebra $\fo_N[t]$ admits a natural structure of {\it restricted Lie algebra}
with the $p$-map $\fo_N[t]\rightarrow \fo_N[t]$ sending $x\mapsto x^{[p]}$.
More precisely, it is a restricted subalgebra of $\gl_N[t]$.
Then the $p$-map is defined on a basis of $\fo_N[t]$ by the role $(F_{i,j}t^r)^{[p]}=\delta_{i,j}F_{i,j}t^{rp}$ (see Section \ref{section current Lie}).
By general theory (cf. \cite[Section 2.3]{Jan98}), the enveloping algebra $\U(\fo_N[t])$ has a large {\it $p$-center} $Z_p(\U(\fo_N[t]))$ generated by the elements
\[
\{(F_{i,j}t^r)^p-\delta_{i,j}F_{i,j}t^{rp};~1\leq i,j\leq N, i+j<N+1, r\geq 0\}.
\]
In Section \ref{section current Lie},
we determine the center of the enveloping algebra $\U(\fo_N[t])$.
It is natural to look for lifts of the $p$-central elements in $Z(\xon)$.
Using the relations established in Section \ref{section:OY},
we investigate the $p$-central elements.
We will give a description of the $p$-center $Z_p(\xon)$ and obtain the precise formulas for the generators.
It is shown in Section \ref{section:center} that the generators of
$Z_p(\xon)$ provide the lifts of generators for $Z_p(\U(\fo_N[t]))$.
With this information in hand, we show
in particular that the center $Z(\xon)$ is generated by  $Z_{\HC}(\xon)$ and $Z_p(\xon)$.
We also determine the center for the Yangian $\yon$.
It would also be interesting to generalize the results in this article to the modular Yangians of type $C$,
but it requires more laborious calculations and some new techniques.

\bigskip

\emph{Throughout this paper, $\kk$ denotes an algebraically closed field of characteristic $\Char(\kk)=:p>3$.}
\bigskip
\section{The orthogonal current Lie algebra and $p$-center}\label{section current Lie}
Define the orthogonal Lie algebras $\fo_N$ with $N=2n+1$ and $N=2n$ (corresponding to types $B$ and $D$, respectively) as subalgebras of $\gl_N$ spanned by all elements $F_{i,j}$
\begin{align}\label{Fij}
F_{i,j}:=E_{i,j}-E_{j',i'},
\end{align}
where the elements $\{E_{i,j};~i,j=1,\dots, N\}$ denote the standard basis of $\gl_N$.
Here we use the notation $i'=N+1-i$.
It is easy to verify that the elements $F_{i,j}$ satisfy the relations
\[
[F_{i,j},F_{k,l}]=\delta_{k,j}F_{i,l}-\delta_{i,l}F_{k,j}-\delta_{k,i'}F_{j',l}+\delta_{l,j'}F_{k,i'}.
\]

The {\it current algebra} is the Lie algebra $\fo_N[t]:=\fo_N\otimes \kk[t]$.
When $x\in\fo_N$ and $f\in\kk[t]$ we usually abbreviate $x\otimes f=xf\in \fo_N[t]$.
As a vector space, $\fo_N[t]$ is spanned by elements $\{F_{i,j}t^r;~1\leq i,j\leq N, r\geq 0\}$,
and the Lie bracket is given by
\begin{align}\label{bracket Fijtr}
[F_{i,j}t^r,F_{k,l}t^s]=(\delta_{k,j}F_{i,l}-\delta_{i,l}F_{k,j}-\delta_{k,i'}F_{j',l}+\delta_{l,j'}F_{k,i'})t^{r+s}.
\end{align}
We observe that $F_{i,j}+F_{j',i'}=0$,
so that the elements $F_{i,j}t^r$ with $r=0,1,2,\dots$ and $i+j<N+1$ make a basis of $\fo_N[t]$.

It is well-known that $\fo_N$ is a restricted simple Lie subalgebra of $\gl_N$ (cf. \cite[Section 0.13]{Hum95}).
Moreover,
the current algebra $\fo_N[t]$ is a restricted Lie algebra with $p$-map defined on the
basis by the rule $(xt^r)^{[p]}:=x^{[p]}t^{rp}$ where $x^{[p]}$ denotes
the $p$th matrix power of $x\in\fo_N$ (see \cite[Lemma 3.3]{BT18}).
In particular, one can show by direct computation that $(F_{i,j}t^r)^{[p]}=\delta_{i,j}F_{i,j}t^{rp}$.

Let $Z(\fo_N[t]):=Z(\U(\fo_N[t]))$ be the center of the universal enveloping algebra $\U(\fo_N[t]))$.
Using the restricted structure,
we can define the {\it $p$-center} $Z_p(\fo_N[t])$ of $\U(\fo_N[t])$ to be the subalgebra of $Z(\fo_N[t])$ generated by $x^p-x^{[p]}$ for all $x\in\fo_N[t]$.
Since the map $x\mapsto x^p-x^{[p]}$ is $p$-semilinear,
we have that
\begin{align}\label{zp-ont}
Z_p(\fo_N[t])=\kk\big[(F_{i,j}t^r)^p-\delta_{i,j}F_{i,j}t^{rp};~1\leq i,j\leq N, i+j<N+1, r\geq 0\big]
\end{align}
as a free polynomial algebra.

There is a filtration
\[
\U(\fo_N[t])=\bigcup\limits_{r\geq  0}{\rm F}_r\U(\fo_N[t])
\]
of the universal enveloping algebra $\U(\fo_N[t])$,
which is defined by by placing $F_{i,j}t^r$ in degree $r+1$.
Then the associated graded algebra $\gr\U(\fo_N[t])$ is isomorphic (both as a graded algebra and as a graded $\fo_N[t]$-module)
to the symmetric algebra $S(\fo_N[t])$.
The adjoint action of $\fo_N[t]$ on itself extends to actions of $\fo_N[t]$ on $\U(\fo_N[t])$
and $S(\fo_N[t])$ by derivations. The invariant subalgebras are denoted by
$\U(\fo_N[t])^{\fo_N[t]}$ and $S(\fo_N[t])^{\fo_N[t]}$, respectively.
It is clear that
\begin{align}\label{grZ}
\gr Z(\fo_N[t])\subseteq S(\fo_N[t])^{\fo_N[t]}
\end{align}

The following result is adapted from \cite[Lemma 3.2]{BT18}.
We provide a proof here which is slightly different in the type $B$ case.
\begin{Lemma}\label{Lemma: center of U}
The invariant algebra $S(\fo_N[t])^{\fo_N[t]}$ is freely generated by
\[
 S(\fo_N[t])^{\fo_N[t]}~\text{is freely generated by}~\{(F_{i,j}t^r)^p;~~1\leq i,j\leq N, i+j<N+1, r\geq 0\}.
\]
\end{Lemma}
\begin{proof}
Note that all the elements in the set $\{(F_{i,j}t^r)^p;~~1\leq i,j\leq N, i+j<N+1, r\geq 0\}$ belong to $S(\fo_N[t])^{\fo_N[t]}$
and let $I(\fo_N[t])$ be the subalgebra of $S(\fo_N[t])^{\fo_N[t]}$ generated by them.
Let
\[
B:=\{(i,j,r);~1\leq i,j\leq N, i+j<N+1, r\geq 0\}
\]
for short.
It follows that
\[
S(\fo_N[t])=\kk[F_{i,j}t^r;~(i,j,r)\in B]~\text{and}~I(\fo_N[t])=\kk[(F_{i,j}t^r)^p;~(i,j,r)\in B].
\]
Hence, $S(\fo_N[t])$ is free as an $I(\fo_N[t])$-module with basis $\{\Pi_{(i,j,r)\in B}(F_{i,j}t^r)^{\omega(i,j,r)};~\omega\in\Omega\}$,
where
\[
\Omega:=
\left\{
\omega:B\rightarrow \NN;~~
\begin{array}{l}
0 \leq \omega(i,j,r)<p,~\forall (i,j,r)\in B,\\
\text{$\omega(i,j,r)=0$ for all but finitely many}~(i,j,r)\in
B
\end{array}
\right\}.
\]
To complete the proof of the lemma,
we must show that $S(\fo_N[t])^{\fo_N[t]}\subseteq I(\fo_N[t])$.
Given $f\in S(\fo_N[t])^{\fo_N[t]}$,
we thus write
\[
f=\sum\limits_{w\in\Omega}c_{\omega}\prod\limits_{(i,j,r)\in B}(F_{i,j}t^r)^{\omega(i,j,r)}
\]
for $c_\omega\in I(\fo_N[t])$, all but finitely many of which are zero.
Also fix a non-zero function $\omega$, we have to prove that $c_\omega=0$.

Suppose first that $\omega(i,j,r)>0$ for some $(i,j,r)\in B$ with $i\neq j$.
Choose $s\in\NN$ that it is bigger that all $h$ such that $\omega(k,l,h)>0$ for $(k,l,h)\in B$.
If $i\neq i'$, then we see that
\begin{align*}
\ad(F_{i,i}t^s)(f) =\sum_{\omega\in\Omega}c_{\omega}\sum_{\substack{(k,l,h)\in B \\ \omega(k,l,h)> 0}}\omega(k,l,h)&(F_{k,l} t^{h})^{\omega(k,l,h)-1}\left[F_{i,i}t^s, F_{k,l}t^{h}\right]\\
&\times\prod_{\substack{(k'',l'',h'') \in B \\ (k'',l'',h'')\neq (k,l,h)}}
(F_{k'',l''}t^{h''})^{\omega(k'',l'',h'')},
\end{align*}
The choice of $s$ in conjunction with $[F_{i,i}t^s, F_{i,j}t^{r}]=F_{i,j}t^{s+r}$ implies that the coefficient of
\[
(F_{i,j}t^r)^{\omega(i,j,r)-1}F_{i,j}t^{r+s}\prod_{\substack{(k'',l'',h'') \in B \\ (k'',l'',h'')\neq (i,j,r)}}
(F_{k'',l''}t^{h''})^{\omega(k'',l'',h'')}
\]
in this expression is $c_\omega\omega(i,j,r)$.
It must be zero because $f\in S(\fo_N[t])^{\fo_N[t]}$.
Since $0<\omega(i,j,r)<p$,
we conclude that $c_\omega=0$.
If $i=i'$, then we consider the coefficient of the same term in $\ad(F_{j,j}t^s)(f)$ for $s$ as before.
By the same token, we can treat the case that $\omega(j,j,r)>0$ for some $(j,j,r)\in B$.
\end{proof}

\begin{Theorem}\label{Thm: center of U}
We have that the center $Z(\fo_N[t])$ of $\U(\fo_N[t])$ is equal to the $p$-center $Z_p(\fo_N[t])$.
In particular, $Z(\fo_N[t])$ is freely generated by
\begin{align}\label{gene of center}
\{(F_{i,j}t^r)^p-\delta_{i,j}F_{i,j}t^{rp};~1\leq i,j\leq N, i+j<N+1, r\geq 0\}.
\end{align}
\end{Theorem}
\begin{proof}
Since the degree of the element $(F_{i,j}t^r)^p-\delta_{i,j}F_{i,j}t^{rp}\in\U(\fo_N[t])$ is $rp+p$,
we denote its image in the graded algebra $\gr\U(\fo_N[t])$ by $\gr_{rp+p}\big((F_{i,j}t^r)^p-\delta_{i,j}F_{i,j}t^{rp}\big)$.
It follows that
\[\gr_{rp+p}\big((F_{i,j}t^r)^p-\delta_{i,j}F_{i,j}t^{rp}\big)=(F_{i,j}t^r)^p\in S(\fo_N[t]).\]
Lemma \ref{Lemma: center of U} implies that the elements \eqref{gene of center} are lifts of the algebraically independent generators of $S(\fo_N[t])^{\fo_N[t]}$.
We obtain $S(\fo_N[t])^{\fo_N[t]}\subseteq\gr Z_p(\fo_N[t])$.
Owing to \eqref{grZ},
we also have $\gr Z_p(\fo_N[t])\subseteq \gr Z(\fo_N[t])\subseteq S(\fo_N[t])^{\fo_N[t]}$,
so equality must hold throughout:
$ \gr Z(\fo_N[t])=\gr Z_p(\fo_N[t])=S(\fo_N[t])^{\fo_N[t]}$.
This yields $ Z(\fo_N[t])=Z_p(\fo_N[t])$.
\end{proof}

\section{Orthogonal Yangians}\label{section:OY}
\subsection{The RTT generators}\label{subsection title-Basic}
Following \cite{AACFR03}, we define the {\it extended Yangian} $\X(\fo_N)$ as an associative algebra over $\kk$ with {\it RTT generators} $\{t_{i,j}^{(r)};~1\leq i,j\leq N, r\geq 1\}$ satisfying certain quadratic relations.
In order to write them down, introduce the formal series
\begin{align}\label{tiju}
t_{i,j}(u):=\delta_{i,j}+\sum\limits_{r\geq 1}t_{i,j}^{(r)}u^{-r}\in\xon[[u^{-1}]]
\end{align}
and combine them into the matrix $T(u):=(t_{i,j}(u))_{1\leq i,j\leq N}$.
In terms of these series, the defining relations for the algebra $\xon$ are written as
\begin{align}\label{def relations RTT}
[t_{i,j}(u),t_{k,l}(v)]&=\frac{1}{u-v}\big(t_{k,j}(u)t_{i,l}(v)-t_{k,j}(v)t_{i,l}(u)\big)\\
&-\frac{1}{u-v-\kappa}\big(\delta_{k,i'}\sum\limits_{q=1}^N t_{q,j}(u)t_{q',l}(v)-\delta_{l,j'}\sum\limits_{q=1}^N t_{k,q'}(v)t_{i,q}(u)\big),\nonumber
\end{align}
where we set $\kappa:=N/2-1$.

It was shown in \cite{AACFR03} (see also \cite{AMR06}) that the product $T(u-\kappa)T^t(u)$ is a scalar matrix with
\begin{align}\label{HC center}
T(u-\kappa)T^t(u)=c(u)1,
\end{align}
where $c(u):=1+\sum_{r\geq 1}c^{(r)}u^{-r}$ is a series in $u^{-1}$ and $T^t(u)$ denotes the matrix transposition defined by $T^t(u)_{i,j}=t_{j',i'}(u)$.
All its coefficients belong to the center $Z(\X(\fo_N))$ of $\xon$.
The algebra generated by the coefficients $\{c^{(r)};~r>0\}$ will be denoted $Z_{\HC}(\xon)$.
We call it the {\it Harish-Chandra center} of $\xon$.

\subsection{The PBW theorem}
Introduce an ascending filtration on the extended Yangian $\xon$ by setting $\deg t_{i,j}^{(r)}=r-1$ for all $r\geq 1$.
Denote by $\gr_{r-1}t_{i,j}^{(r)}$ the images of the elements $t_{i,j}^{(r)}$ in the $(r-1)$-th component of the associated graded algebra $\gr\xon$.
By the PBW theorem (\cite[Corollary 3.10]{AMR06}),
the mapping
\begin{align}\label{PBW map}
\gr_{r-1}t_{i,j}^{(r)}\mapsto F_{i,j}t^{r-1}+\frac{1}{2}\delta_{i,j}\zeta_r
\end{align}
defines an isomorphism
\[
\gr\xon\cong\U(\fo_N[t])\otimes \kk[\zeta_1,\zeta_2,\dots],
\]
where $\kk[\zeta_1,\zeta_2,\dots]$ is the algebra of polynomials in indeterminants $\zeta_r$ and
$\zeta_r$ is the image of $\gr_{r-1}c^{(r)}$.

\subsection{Gaussian generators}
Apply the Gauss deomposition to the generator matrix $T(u)$ associated with the extended Yangian $\xon$,
\begin{align}\label{Gauss decomp}
T(u)=F(u)H(u)E(u),
\end{align}
where $F(u), H(u)$ and $E(u)$ are uniquely determined matrices of the form
$$
H(u) = \left(
\begin{array}{cccc}
h_{1}(u) & 0&\cdots&0\\
0 & h_{2}(u) &\cdots&0\\
\vdots&\vdots&\ddots&\vdots\\
0&0 &\cdots&h_{N}(u)
\end{array}
\right),
$$$$
E(u) =
\left(
\begin{array}{cccc}
1 & e_{1,2}(u) &\cdots&e_{1,N}(u)\\
0 & 1 &\cdots&e_{2,N}(u)\\
\vdots&\vdots&\ddots&\vdots\\
0&0 &\cdots&1
\end{array}
\right),\:
F(u) = \left(
\begin{array}{cccc}
1 & 0 &\cdots&0\\
f_{2,1}(u) & 1 &\cdots&0\\
\vdots&\vdots&\ddots&\vdots\\
f_{N,1}(u)&f_{N,2}(u) &\cdots&1
\end{array}
\right).
$$
In terms of quasideterminants of \cite{GR97} (see also \cite[Section 1.11]{Mol07}),
we have the following descriptions (cf. \cite[Section 4]{JLM18}):
\begin{align}\label{quasideterminants H}
h_i(u) =
\left|
\begin{array}{cccc}
t_{1,1}(u) & \cdots & t_{1,i-1}(u)&t_{1,i}(u)\\
\vdots & \ddots &\vdots&\vdots\\
t_{i-1,1}(u)&\cdots&t_{i-1,i-1}(u)&t_{i-1,i}(u)\\
t_{i,1}(u) & \cdots & t_{i,i-1}(u)&
\hbox{\begin{tabular}{|c|}\hline$t_{i,i}(u)$\\\hline\end{tabular}}
\end{array}
\right|,
\end{align}
whereas
\begin{equation}\label{quasideterminants E}
e_{i,j}(u) =
h_i(u)^{-1} \left|
\begin{array}{cccc}
t_{1,1}(u) & \cdots &t_{1,i-1}(u)& t_{1,j}(u)\\
\vdots & \ddots &\vdots&\vdots\\
t_{i-1,1}(u) & \cdots &t_{i-1,i-1}(u)&t_{i-1,j}(u)\\
t_{i,1}(u) & \cdots & t_{i,i-1}(u)&
\hbox{\begin{tabular}{|c|}\hline$t_{i,j}(u)$\\\hline\end{tabular}}
\end{array}
\right|,
\end{equation}
and
\begin{equation}\label{quasideterminants F}
f_{j,i}(u) =
\left|
\begin{array}{cccc}
t_{1,1}(u) & \cdots &t_{1,i-1}(u)& t_{1,i}(u)\\
\vdots & \ddots &\vdots&\vdots\\
t_{i-1,1}(u) & \cdots & t_{i-1,i-1}(u)&t_{i-1,i}(u)\\
t_{j,1}(u) & \cdots & t_{j,i-1}(u)&
\hbox{\begin{tabular}{|c|}\hline$t_{j,i}(u)$\\\hline\end{tabular}}
\end{array}
\right|{h}_i(u)^{-1}.
\end{equation}
We use the following notation for the coefficients:
\[
h_i(u)=\sum\limits_{r\geq 0}h_i^{(r)}u^{-r},~\tilde{h}_i(u)=\sum\limits_{r\geq 0}\tilde{h}_i^{(r)}u^{-r}:=h_i(u)^{-1};
\]
\[
e_{i,j}(u)=\sum\limits_{r\geq 1}e_{i,j}^{(r)}u^{-r},~f_{j,i}(u)=\sum\limits_{r\geq 1}f_{j,i}^{(r)}u^{-r}.
\]
Furthermore, set
\begin{align}\label{def:kef-1}
k_i(u):=\tilde{h}_i(u)h_{i+1}(u),~e_i(u):=e_{i,i+1}(u),~f_{i}(u):=f_{i+1,i}(u),
\end{align}
with $i=1,\dots,n$ for type $B$ and with $i=1,\dots,n-1$ for type $D$.
In the latter case we also set
\begin{align}\label{def:kn for D}
k_n(u):= \tilde{h}_{n-1}(u)h_{n+1}(u)
\end{align}
and
\begin{align}\label{def:efn for D}
e_{n}(u):=e_{n-1,n+1}(u),~f_n(u):=f_{n+1,n-1}(u).
\end{align}
We will also use the coefficients of the series defined by
\begin{align}\label{coeff of eifi}
e_{i}(u)=\sum\limits_{r\geq 1}e_{i}^{(r)}u^{-r}~\text{and}~f_{i}(u)=\sum\limits_{r\geq 1}f_{i}^{(r)}u^{-r}.
\end{align}

\subsection{Drinfeld presentation of the extended Yangian}
We will give the modular analogue of \cite[Theorem 5.14]{JLM18}.
To state the result, we will assume that the simple roots of $\fo_N$ are $\alpha_1,\dots,\alpha_n$ with $\alpha_i:=\epsilon_i-\epsilon_{i+1}$ for $i=1,\dots,n-1$,
and
\[
\alpha_n:=\left\{
\begin{array}{ll}
\epsilon_n &~\text{for}~N=2n+1,\\
\epsilon_{n-1}+\epsilon_n & ~\text{for}~N=2n,
\end{array}{}
\right.
\]
where $\epsilon_1,\dots,\epsilon_n$ is an orthogonal basis of a vector space with the bilinear form such that $(\epsilon_i,\epsilon_i)=1$ for $i=1,\dots,n$.
The associated Cartan matrix $C=(c_{i,j})_{i,j=1}^n$ is defined by $c_{i,j}:=(\alpha_i,\alpha_j)$ for series $D$, and by
\[
c_{i,j}:=\left\{
\begin{array}{ll}
(\alpha_i,\alpha_j) &~\text{if}~i<n,\\
2(\alpha_i,\alpha_j) & ~\text{if}~i=n,
\end{array}{}
\right.
\]
for series $B$.
We will use the series introduced in \eqref{def:kef-1}-\eqref{def:efn for D} along with
\[
e^{\circ}_i(u):=\sum\limits_{r>1}e_i^{(r)}u^{-r}\ \ \ \text{and} \ \ \ f^{\circ}_i(u):=\sum\limits_{r>1}f_i^{(r)}u^{-r}.
\]
\begin{Theorem}\label{Theorem: Drinfeld presentation}
The extended Yangian $\X(\fo_N)$ is generated by the coefficients of the series $h_i(u)$ with $i=1,\dots,n+1$, and the series
$e_i(u)$ and $f_i(u)$ with $i=1,\dots, n$,
subject only to the following relations, where the indices
take all admissible values unless specified otherwise.
We have

\begin{align} \big[h_i(u),h_j(v)\big]=&\,0,\label{hihj}\\
\big[e_i(u),f_j(v)\big]=&\,\delta_{i,j}\frac{k_i(u)-k_i(v)}{u-v}.\label{eifj}
\end{align}
For $i\leq n$ and all $j$,
and for $i=n+1$ and $j\leq n-2$ we have
\begin{align}
\big[h_i(u),e_j(v)\big]=&-(\epsilon_i,\alpha_j)\frac{h_i(u)\big(e_j(u)-e_j(v)\big)}{u-v},\label{hiej}\\
\big[h_i(u),f_j(v)\big]=&\,(\epsilon_i,\alpha_j)\frac{\big(f_j(u)-f_j(v)\big)h_i(u)}{u-v},\label{hifj}
\end{align}
where $\epsilon_{n+1}=0$.
For $N=2n+1$, we have
\begin{align}\label{hn+1efn-1}
[h_{n+1}(u),e_{n-1}(v)]=0=[h_{n+1}(u),f_{n-1}(v)],
\end{align}
\begin{align}\label{hn+1en}
\big[h_{n+1}(u),e_{n}(v)\big]&=\frac{1}{2(u-v)}h_{n+1}(u)\big(e_n(u)-e_n(v)\big)\\
{}&-\frac{1}{2(u-v-1)}\big(e_n(u-1)-e_n(v)\big)h_{n+1}(u)\nonumber
\end{align}
and
\begin{align}\label{hn+1fn}
\big[h_{n+1}(u),f_{n}(v)\big]&=-\frac{1}{2(u-v)}\big(f_n(u)-f_n(v)\big)h_{n+1}(u)\\
{}&+\frac{1}{2(u-v-1)}h_{n+1}(u)\big(f_n(u-1)-f_n(v)\big),\nonumber
\end{align}
whereas for $N=2n$ we have
\begin{align}\label{hn+1efn-1-2n}
\big[h_{n+1}(u),e_{n-1}(v)\big]&=\frac{h_{n+1}(u)\big(e_{n-1}(v)-e_{n-1}(u)\big)}{u-v},\\
\big[h_{n+1}(u),f_{n-1}(v)\big]&=-\frac{\big(f_{n-1}(v)-f_{n-1}(u)\big)h_{n+1}(u)}{u-v}\nonumber
\end{align}
and
\begin{align}\label{hn+1efn-2n}
\big[h_{n+1}(u),e_{n}(v)\big]&=\frac{h_{n+1}(u)\big(e_{n}(u)-e_{n}(v)\big)}{u-v},\\
\big[h_{n+1}(u),f_{n}(v)\big]&=-\frac{\big(f_{n}(u)-f_{n}(v)\big)h_{n+1}(u)}{u-v}.\nonumber
\end{align}
Moreover,
\begin{align}
\label{eiei}
&\big[e_i(u),e_{i}(v)\big]=\frac{(\alpha_i,\alpha_{i})}{2}
\frac{\big(e_{i}(u)-e_{i}(v)\big)^2}{u-v},\\
\label{fifi}
&\big[f_i(u),f_{i}(v)\big]=-\frac{(\alpha_i,\alpha_{i})}{2}\frac{\big(f_{i}(u)-f_{i}(v)\big)^2}{u-v},
\end{align}
and for $i<j$ we have
\begin{align}\label{eiej}
u\big[e^{\circ}_i(u),e_{j}(v)\big]-v\big[e_i(u),e^{\circ}_{j}(v)\big]
&=-(\alpha_i,\alpha_{j}) e_{i}(u)e_{j}(v),\\
\label{fifj}
u\big[f^{\circ}_i(u),f_{j}(v)\big]-v\big[f_i(u),f^{\circ}_{j}(v)\big]
&=(\alpha_i,\alpha_{j}) f_{j}(v) f_{i}(u).
\end{align}
We have the Serre relations
\begin{align}
\sum_{\sigma\in\fS_k}\big[e_{i}(u_{\sigma(1)}),
\big[e_{i}(u_{\sigma(2)}),\dots,\big[e_{i}(u_{\sigma(k)}),e_{j}(v)\big]\dots\big]\big]&=0,\label{Serre-e}\\
\sum_{\sigma\in\fS_k}\big[f_{i}(u_{\sigma(1)}),
\big[f_{i}(u_{\sigma(2)}),\dots,\big[f_{i}(u_{\sigma(k)}),f_{j}(v)\big]\dots\big]\big]&=0,\label{Serre-f}
\end{align}
for $i\ne j$  with
$k=1+|c_{ij}|$.
\end{Theorem}
\begin{proof}
Apart from the Serre relations where $i$ or $j$ takes the value $n$, all the relations are satisfied in the Yangian $\xon$ due to \cite[Theorem 5.14]{JLM18},
see also \cite[Theorem 6.1]{Mol24},
the arguments there works perfectly in
positive characteristic.
The proof of the Serre relations is long and
will be given in Appendix \ref{section name:appendix} (see Proposition \ref{prop:app-serre relation}).
\end{proof}
\begin{Remark}
(1) We note that the relations (5.47), (5.56) and (5.57) in \cite[Theorem 5.14]{JLM18} contain some typos,
as corrected in \cite[Page 487]{Mol24}.
The counterpart of \eqref{hn+1fn} therein should be corrected by swapping the factors, while the correct condition
for the counterparts of \eqref{eiej} and \eqref{fifj} should be $i<j$.

(2) We identify $\gr\xon $ and $\U(\fo_N[t])\otimes \kk[\zeta_1,\zeta_2,\dots]$ via the isomorphism \eqref{PBW map}.
The proof in \cite[Theorem 5.14]{JLM18} shows moreover that the images of the elements $e_{i,j}^{(r+1)}, f_{j,i}^{(r+1)}$ and $h_i^{(r+1)}$ all belong to
the $r$-th component of $\gr\xon$,
and under our identification we have that
\begin{equation}\label{identification-1}
F_{i,j}t^{r}+\frac{1}{2}\delta_{i,j}\zeta_{r+1}
=\left\{
\begin{array}{lll}
\gr_{r}e_{i,j}^{(r+1)}&\text{if}&i<j<i',\\
\gr_{r}f_{i,j}^{(r+1)}&\text{if}&j<i<j',\\
\gr_{r}h_{i}^{(r+1)}&\text{if}&1\leq i=j\leq n+1.
\end{array}\right.
\end{equation}
Recall that $F_{i,i}+F_{i',i'}=0$.
Then \eqref{identification-1} readily implies that\footnote{The corrects the formula for type $D$ in \cite[Section 5.5]{JLM18},
see also the footnote in \cite[Page 489]{Mol24}.}
\begin{equation}\label{identification-2}
\gr_r h_{n+1}^{(r+1)}
=\left\{
\begin{array}{lll}
\frac{1}{2}\zeta_{r+1}&\text{for}&N=2n+1,\\
-F_{n,n}t^r+\frac{1}{2}\zeta_{r+1}&\text{for}&N=2n.
\end{array}\right.
\end{equation}
\end{Remark}

\begin{Lemma}
For any $i\neq j$, we have the relation  in the extended Yangian $\xon$:
\begin{align}\label{serre-ef30}
[e_i(u),[e_i(u),[e_i(u),e_j(v)]]]=0=[f_i(u),[[f_i(u),[f_i(u),f_j(v)]]].
\end{align}
\end{Lemma}
\begin{proof}
This is a direct consequence of the Serre relations (\eqref{Serre-e} and \eqref{Serre-f}).
\end{proof}


\begin{Lemma}\label{various formula}
The following relations hold in the extended Yangian $\xon$:
\begin{align}
(u-v)[e_i(u),h_i(v)]&=-h_i(v)\big(e_i(v)-e_i(u)\big),\label{eihi-1}\\
(u-v)[e_i(u),h_{i+1}(v)]&=h_{i+1}(v)\big(e_i(v)-e_i(u)\big),\,\,\,\,\,\,\,\,\,\,\,where~i\neq n~for~N=2n+1,\label{eihi+1-1}\\
(u-v)[e_n(u),h_{n-1}(v)]&=
\left\{
\begin{array}{ll}
0,&for~N=2n+1,\\
-h_{n-1}(v)(e_n(v)-e_n(u)),&for~N=2n,
\end{array}
\right.\label{enhn-1-1}\\
(u-v)[e_i(u),\tilde{h}_i(v)]&=\big(e_i(v)-e_i(u)\big)\tilde{h}_i(v),\label{eihitilde}\\
(u-v)[e_i(u),\tilde{h}_{i+1}(v)]&=-\big(e_i(v)-e_i(u)\big)\tilde{h}_{i+1}(v),\,\,\,\,\,\,where~i\neq n ~for~N=2n+1,\label{eihi+1tilde}\\
(u-v)[e_n(u),\tilde{h}_{n-1}(v)]&=
\left\{
\begin{array}{ll}
0,&for~N=2n+1,\\
(e_n(v)-e_n(u))\tilde{h}_{n-1}(v),&for~N=2n.
\end{array}
\right.\label{enhn-1tilde}
\end{align}
\end{Lemma}
\begin{proof}
Equations \eqref{eihi-1}-\eqref{enhn-1-1} follow immediately from \eqref{hiej} and \eqref{hn+1efn-2n}.
Then \eqref{eihitilde}-\eqref{enhn-1tilde} follow from \eqref{eihi-1}-\eqref{enhn-1-1} using $h_j(v)\tilde{h}_j(v)=1$.
\end{proof}
\begin{Corollary}
The following relations hold in the extended Yangian $\xon$:
\begin{align}
e_i(u-(\epsilon_i,\alpha_i))h_i(u)&=h_i(u)e_i(u),& \  \tilde{h}_i(u)e_i(u-(\epsilon_i,\alpha_i))&=e_i(u)\tilde{h}_i(u),\label{eihihiei-1}\\
e_i(u-(\epsilon_{i+1},\alpha_i))h_{i+1}(u)&=h_{i+1}(u)e_i(u),& \  \tilde{h}_{i+1}(u)e_i(u-(\epsilon_i,\alpha_i))&=e_i(u)\tilde{h}_{i+1}(u),\label{eihihiei-2}
\end{align}
where $i\neq n$ in \eqref{eihihiei-2}.
For $N=2n$ we have
\begin{align}
e_n(u+1)h_{n+1}(u)&=h_{n+1}(u)e_n(u),&\  e_{n-1}(u-1)h_{n+1}(u)&=h_{n+1}(u)e_{n-1}(u),\label{comm in N=2n}
\end{align}
\end{Corollary}
\begin{proof}
Equations \eqref{eihihiei-1}-\eqref{eihihiei-2} follow from \eqref{eihi-1}-\eqref{eihi+1-1} and \eqref{eihitilde}-\eqref{eihi+1tilde} by specializing $v$.
For example, to get the first relation in \eqref{eihihiei-1}, set $v:=u+(\epsilon_i,\alpha_i)$ in \eqref{eihi-1}, simplify,
then replace $u$ by $u-(\epsilon_i,\alpha_i)$.
The relations \eqref{comm in N=2n} follow from \eqref{hn+1efn-1-2n} and \eqref{hn+1efn-2n}.
\end{proof}
\begin{Lemma}\label{useful lemma 1}
The following relations hold in the extended Yangian $\X(\fo_N)$ for all $m \geq 0$:
\begin{align}
(u-v)[e_i(u),\big(e_i(v)\!-\!e_i(u)\big)^m]&=\frac{(\alpha_i,\alpha_i)}{2}m\big(e_i(v)\!-\!e_i(u)\big)^{m+1},\label{new1}\\
(u-v)[e_i(u), h_i(v)\big(e_i(v)\!-\!e_i(u)\big)^m] &=\left(\frac{(\alpha_i,\alpha_i)}{2}m-1\right)h_i(v)\big(e_i(v)\!-\!e_i(u)\big)^{m+1},\label{new2}
\end{align}
\begin{multline}\label{new22}
(u-v)[e_{n-1}(u),h_{n+1}(v)\big(e_{n-1}(v)-e_{n-1}(u)\big)^m]=(m-1)h_{n+1}(v)\big(e_{n-1}(v)-e_{n-1}(u)\big)^{m+1}, \\
for ~N=2n,
\end{multline}
\begin{multline}\label{new222}
(u-v)[e_{n}(u),h_{n-1}(v)\big(e_{n}(v)-e_{n}(u)\big)^m]=(m-1)h_{n-1}(v)\big(e_{n}(v)-e_{n}(u)\big)^{m+1}, \\
for~N=2n,
\end{multline}
\begin{multline}\label{new3}
(u-v)[e_i(u), h_{i+1}(v)(e_i(v)\!-\!e_i(u))^m]
=(m+1)h_{i+1}(v)(e_i(v)\!-\!e_i(u))^{m+1},\\
\mbox{~where~}i\ne n\mbox{~for~} N=2n+1,\end{multline}
\begin{multline}\label{new4-2}
[e_n(u), h_{n+1}(v)\big(e_n(v)\!-\!e_n(u)\big)^m]
=(m+1)\frac{1}{2(u-v)}h_{n+1}(v)\big(e_n(v)\!-\!e_n(u)\big)^{m+1}\\
+\frac{1}{2(v-u-1)}\big(e_n(v-1)\!-\!e_n(u)\big)h_{n+1}(v)\big(e_n(v)\!-\!e_n(u)\big)^{m} \quad\mbox{for~}N=2n+1,
\end{multline}
\begin{multline}\label{new4-3}
[e_n(u), \big(e_n(v-1)\!-\!e_n(u)\big)h_{n+1}(v)\big(e_n(v)\!-\!e_n(u)\big)^m]
=(m+1)\frac{1}{2(u-v)}\times\\
\big(e_n(v-1)\!-\!e_n(u)\big)h_{n+1}(v)\big(e_n(v)\!-\!e_n(u)\big)^{m+1} \quad\mbox{for~}N=2n+1,\end{multline}
\begin{multline}\label{new5}
(u-v)[e_i(u), \,h_{i+1}(v)\big(e_i(v)\!-\!e_i(u)\big)^m\tilde{h}_i(v)]
=(m+2)h_{i+1}(v)\big(e_i(v)\!-\!e_i(u)\big)^{m+1}\tilde{h}_i(v), \\
\mbox{~where~}i\ne n\mbox{~for~} N=2n+1,
\end{multline}
\begin{multline}\label{new6}
[e_n(u), \,h_{n+1}(v)\big(e_n(v)\!-\!e_n(u)\big)^m\tilde{h}_n(v)]
=\frac{(m+3)}{2(u-v)}h_{n+1}(v)\big(e_n(v)\!-\!e_n(u)\big)^{m+1}\tilde{h}_n(v)\\+\frac{1}{2(v-u-1)}\big(e_n(v-1)\!-\!e_n(u)\big)h_{n+1}(v)\big(e_n(v)\!-\!e_n(u)\big)^{m}\tilde{h}_n(v)\quad\mbox{for~}N=2n+1,
\end{multline}
\begin{multline}\label{new7}
[e_n(u), \,h_{n+1}(v)\big(e_n(v)\!-\!e_n(u)\big)^m\tilde{h}_{n-1}(v)]
=(m+2)h_{n+1}(v)\big(e_n(v)\!-\!e_n(u)\big)^{m+1}\tilde{h}_{n-1}(v),\\
for ~N=2n.
\end{multline}
\end{Lemma}
\begin{proof}
The relation \eqref{new1} follows from \eqref{eiei} and the Leibniz rule.
Then \eqref{eihi-1}, \eqref{hn+1efn-1-2n}, \eqref{enhn-1-1}, \eqref{eihi+1-1} and \eqref{hn+1en} imply that \eqref{new2}-\eqref{new4-2} hold for the case $m=0$.
In the general case, we use \eqref{new1} together with the Leibniz rule.
For \eqref{new4-3}, we have by \eqref{new1} and \eqref{hn+1en} that
\begin{align*}
&[e_n(u), \big(e_n(v-1)\!-\!e_n(u)\big)h_{n+1}(v)]\nonumber\\
=&[e_n(u), \big(e_n(v-1)\!-\!e_n(u)\big)]h_{n+1}(v)+\big(e_n(v-1)\!-\!e_n(u)\big)[e_n(u), h_{n+1}(v)]\nonumber\\
=&\frac{1}{2(u-v+1)}\big(e_n(v-1)\!-\!e_n(u)\big)^2h_{n+1}(v)+\big(e_n(v-1)\!-\!e_n(u)\big)\times\nonumber\\
&\left(\frac{1}{2(u-v)}h_{n+1}(v)\big(e_n(v)-e_n(u)\big)+\frac{1}{2(v-u-1)}\big(e_n(v-1)-e_n(u)\big)h_{n+1}(v)\right)\nonumber\\
=&\frac{1}{2(u-v)}\big(e_n(v-1)\!-\!e_n(u)\big)h_{n+1}(v)\big(e_n(v)-e_n(u)\big).
\end{align*}
Then the general case follows from \eqref{new1} and \eqref{new4-2} using Leibniz again.
The relations \eqref{new5}-\eqref{new7} can be treated similarly.
\end{proof}
We introduce the following notation:
\begin{align*}
\mathcal{H}^m_{v,u}:=&\frac{m!}{\big(2(u-v)\big)^m}h_{n+1}(v)\big(e_n(v)\!-\!e_n(u)\big)^m\\
&+\frac{m(m-1)!}{\big(2(u-v)\big)^{m-1}\times 2(v-u-1)}\big(e_n(v-1)\!-\!e_n(u)\big)h_{n+1}(v)\big(e_n(v)\!-\!e_n(u)\big)^{m-1},\\
\tilde{\mathcal{H}}^{m-1}_{v,u}:=&\frac{(m+1)!}{2\times \big(2(u-v)\big)^{m-1}}h_{n+1}(v)\big(e_n(v)\!-\!e_n(u)\big)^{m-1}\tilde{h}_n(v)\\
&+\frac{(m-1)m!}{2\times\big(2(u-v)\big)^{m-2}\times 2(v-u-1)}\big(e_n(v-1)\!-\!e_n(u)\big)h_{n+1}(v)\big(e_n(v)\!-\!e_n(u)\big)^{m-2}\tilde{h}_n(v)
\end{align*}

Using \eqref{new4-2}, \eqref{new4-3} and \eqref{new6},
we obtain the following corollary:
\begin{Corollary}
For $N=2n+1$ we have
\begin{align}\label{en Hmvu}
[e_n(u),\mathcal{H}^m_{v,u}]=\mathcal{H}^{m+1}_{v,u},\ \ \  [e_n(u),\tilde{\mathcal{H}}^{m-1}_{v,u}]=\tilde{\mathcal{H}}^{m}_{v,u}.
\end{align}
\end{Corollary}

We introduce one more family of elements.
Set

\begin{align}\label{definition:down formula}
h_{i\downarrow \ell}(u)&:=h_i(u)h_i(u-1)\cdots h_i(u-\ell+1),
\end{align}
\begin{align}\label{definition:up formula}
h_{i\uparrow \ell}(u)&:=h_i(u)h_i(u+1)\cdots h_i(u+\ell-1).
\end{align}
\begin{Lemma}
The following relations hold for all $i\leq n$ and $\ell\geq 1$:
\begin{align}
(u-v)[h_{i\downarrow \ell}(u),e_i(v)]&=\ell h_{i\downarrow \ell}(u)\big( e_i(v)-e_i(u)\big),\label{down 1}\\
(u-v)[h_{i\uparrow \ell}(u),e_{i-1}(v)]&=\ell h_{i\uparrow \ell}(u)\big( e_{i-1}(u)-e_{i-1}(v)\big).\label{up 1}
\end{align}
For $N=2n$ we have
\begin{align}
(u-v)[h_{n+1\uparrow \ell}(u),e_n(v)]&=\ell h_{n+1\uparrow \ell}(u)\big( e_n(u)-e_n(v)\big),\label{up 2-2n}\\
(u-v)[h_{n+1\downarrow \ell}(u),e_{n-1}(v)]&=\ell h_{n+1\downarrow \ell}(u)\big( e_{n-1}(v)-e_{n-1}(u)\big).\label{down 2-2n}
\end{align}
\end{Lemma}
\begin{proof}
The proof is similar to the proof of \cite[Lemma 4.10]{BT18} (see also \cite[Lemma 3.15]{CH23}).
We prove \eqref{up 2-2n}  in detail here in order to highlight the minor differences.
We will prove it in the following equivalent form:
\begin{align}\label{equiv to up2-2n}
(u-v+\ell)h_{n+1\uparrow \ell}(u)e_n(v)=(u-v)e_n(v)h_{n+1\uparrow \ell}(u)+\ell h_{n+1\uparrow \ell}(u)e_n(u).
\end{align}
This follows when $\ell=1$ from \eqref{hn+1efn-2n}.
To prove \eqref{equiv to up2-2n} in general, we proceed by induction on $\ell$.
Given \eqref{equiv to up2-2n} for some $\ell\geq 1$,
multiply both sides on the left by $(u-v+\ell+1)h_{n+1}(u+\ell)$ to deduce that:
\begin{align}\label{equiv to up2-2n-1}
(u-v+\ell+1)(u-v+\ell)h_{n+1\uparrow \ell+1}(u)e_n(v)&=(u-v)(u-v+\ell+1)h_{n+1}(u+\ell)e_n(v)h_{n+1\uparrow \ell}(u)\nonumber\\
&+\ell(u-v+\ell+1)h_{n+1\uparrow \ell+1}(u)e_n(u).
\end{align}
Using the case of $\ell=1$ in \eqref{equiv to up2-2n} and replacing $u$ by $u+\ell$ give that
\[
(u-v+\ell+1)h_{n+1}(u+\ell)e_n(v)=(u-v+\ell)e_n(v)h_{n+1}(u+\ell)+h_{n+1}(u+\ell)e_n(u+\ell).
\]
Then substituting the above identity into \eqref{equiv to up2-2n-1} and using the first relation in \eqref{comm in N=2n}
we obtain \eqref{equiv to up2-2n} with $\ell$ replaced by
$\ell+1$, as required.
By using the second relation in \eqref{comm in N=2n}
instead of the first one, \eqref{down 2-2n} is a similar argument to \eqref{up 2-2n}.
\end{proof}
\subsection{Automorphisms}\label{subsec: anto}
We list the following (anti)automorphisms of $\xon$ which are needed in the next section.
\begin{enumerate}
\item(``Transposition'')
The mapping $\tau: t_{i,j}(u)\mapsto t_{j,i}(u)$ defines an anti-automorphism of $\xon$.
On Drinfeld generators, we have
\[\tau(e_{i,j}(u))=f_{j,i}(u), \tau(f_{j,i}(u))=e_{i,j}(u)\]
for $i<j$, and $\tau(h_i(u))=h_i(u)$ for all $i$ (cf. \cite[(2.9) and Lemma 4.1]{Mol24}).
\item(``Permutation'') Let $S_{N}$ be the Symmetric group on $N$ objects.
For each $w\in S_N$, if $w(i')=w(i)'$ for all $1\leq i\leq N$,
then there is an automorphism
$w:\xon\rightarrow \xon$ sending $t_{i,j}^{(r)} \mapsto t_{w(i),
  w(j)}^{(r)}$. This is clear from the RTT relation (\ref{def relations RTT}).
\end{enumerate}

For later use we define
\[
I_1:=\{(i,j);~1\leq i,j\leq N,~i<j<i'\}\backslash\{(i,n+1);~1\leq i<n\},\ \ I_2:=\{(n+1,j);~n'<j\leq 1'\}.
\]
for $N=2n+1$, and we set $I:=I_1\cup I_2$.
For $N=2n$ we set
\[
I:=\{(i,j);~1\leq i,j\leq N,~i<j<i'\}.
\]
Moreover, we also set $J:=\{(j,i);~1\leq i,j\leq N, (i,j)\in I\}$.

\begin{Lemma}\label{permutation auto}
If $(i,j)\in I_1$ or $(i,j)\in I$ for $N=2n$,
then there exists a permutation automorphism $\omega$ of $\xon$ sending
$e_{i,j}(u)\mapsto e_i(u)$.
If $(i,j)\in I_2$, then there is a permutation automorphism $\omega$ sending $e_{i,j}(u)\mapsto e_{n+1,n+2}(u)$.
\end{Lemma}
\begin{proof}
We first assume that $N=2n+1$ and $(i,j)\in I_1$.
By definition, we may assume further $1\leq i\leq n-1$.
As $j\neq n+1$, the transposition $(j,i+1)$ can be extended to an automorphism of $\xon$.
Then \eqref{quasideterminants H}-\eqref{quasideterminants F} implies the transposition $(j,i+1)$ maps $e_{i,j}(u)\mapsto e_i(u)$.
Let $(i,j)\in I_2$.
We then extend the transposition $(j,n+2)$ to a permutation $\omega=(j,n+2)(j',n)$.
It is easy to see that $\omega$ gives rise to a permutation automorphism sending $e_{n+1,j}(u)\mapsto e_{n+1,n+2}(u)$ using again \eqref{quasideterminants H}-\eqref{quasideterminants F}.
One argues similarly for the case $N=2n$.
\end{proof}

\begin{Remark}
In fact, the proof of Theorem \ref{Theorem: Drinfeld presentation} also works for the type $C$ Yangians.
However,
to deal with symplectic Lie algebras,
we know that $F_{i,i'}:=2E_{i,i'}\neq 0$ (see for example \cite[(2.1)]{JLM18}).
Hence, we need to consider the Drinfeld generators $e_{i,i'}(u)$.
It seems to be hard to find a suitable automorphism which sends the $e_{i,i'}(u)$ to some $e_j(u)$.
\end{Remark}

\section{The centers of $\xon$ and $\yon$}\label{section:center}
In this section, we will describe the center of the extended Yangian $\xon$,
and give precise formulas for the generators.
\subsection{Harish-Chandra center}
Recall from Subsection \ref{subsection title-Basic} the {\it Harish-Chandra center}
$Z_{\HC}(\xon)$ of the extended Yangian $\xon$
is defined to be the subalgebra generated by the coefficients of $c(u)$ (see \eqref{HC center}).
The following result provides a formula for the series $c(u)$ in terms of the Gaussian generators $h_i(u)$
with $i=1,\dots,n+1$ (\cite[Theorem 5.8]{JLM18} and \cite[Theorem 5.3]{Mol24}).
\begin{Theorem}\label{thm: HC center via Drinfeld}
We have the relations in the extended Yangian $\xon$:
\begin{align*}
c(u)&=\prod\limits_{i=1}^n\frac{h_i(u-i+1)}{h_i(u-i)}\times h_{n+1}(u-n+1/2)h_{n+1}(u-n) \ \ \ \ for ~N=2n+1, and \\
c(u)&=\prod\limits_{i=1}^{n-1}\frac{h_i(u-i+1)}{h_i(u-i)}\times h_{n}(u-n+1)h_{n+1}(u-n+1) \ \ \ \ for ~N=2n.
\end{align*}
\end{Theorem}
\subsection{Off-diagonal $p$-central elements}\label{subs: off-dia}
In this subsection,
we investigate the $p$-central elements that lie in the ``root subalgebras'' $\xon_{i,j}^+, \xon_{j,i}^-$ for $(i,j)\in I$,
that is, the subalgebras generated by $\{e_{i,j}^{(r)};~r>0\}$ and $\{f_{j,i}^{(r)};~r>0\}$, respectively.

\begin{Lemma}\label{offdig-power}
For $(i,j)\in I$, all coefficients in the power series $(e_{i,j}(u))^p$ and $(f_{j,i}(u))^p$ belong to $Z(\xon)$.
\end{Lemma}
\begin{proof}
In view of \cite[Proposition 5.7]{JLM18},
we know that $e_{n+1,n+2}(u)=-e_n(u-1/2)$ for $N=2n+1$.
Then using Lemma \ref{permutation auto} and the anti-automorphism $\tau$,
it only needs to be proved that the coefficients of $(e_i(u))^{p}$ are central in $\xon$ for each $i=1,\dots,n$.
Since  we are in characteristic $p$, it suffices to establish the following identities in $X(\fo_N)[[u^{-1},v^{-1}]]$ for all admissible $j$:
\begin{align}
(\ad e_i(u))^p(e_j(v))&=0,  \label{center-ee} \\
(\ad e_i(u))^p(h_j(v))&=0,  \label{center-eh}\\
(\ad e_i(u))^p(f_j(v))&=0.  \label{center-ef}
\end{align}

When $i\neq j$, the identity \eqref{center-ee} is clear because $(\ad e_i(u))^3(e_j(v))=0$ by \eqref{serre-ef30}.
To show that $(\ad e_i(u))^p(e_i(v))=0$,
we use \eqref{eiei} and \eqref{new1} repeatedly:
\begin{align*}
(u-v)^p(\ad e_i(u))^p(e_i(v))=&\frac{(\alpha_i,\alpha_i)}{2}(u-v)^{p-1}(\ad e_i(u))^{p-1}(e_i(v))\!-\!(e_i(u))^2\\
=&2(\frac{(\alpha_i,\alpha_i)}{2})^2(u-v)^{p-2}(\ad e_i(u))^{p-2}(e_i(v))\!-\!(e_i(u))^3\\
=&\cdots=p!(\frac{(\alpha_i,\alpha_i)}{2})^p(e_i(v)-e_i(u))^{p+1}=0.
\end{align*}

To check \eqref{center-eh}, we first consider the case that $i\leq n-2$.
It is immediate from \eqref{hiej} if $j<i$ or $j>i+1$.
For $j=i$, we have by \eqref{hiej} and \eqref{new2} with $m=1$ that $(\ad e_i(u))^2(h_i(v))=0$.
Then we use \eqref{new3} to see that
\begin{align*}
(u-v)^p(\ad e_i(u))^p(h_{i+1}(v))&=(u-v)^{p-1}(\ad e_i(u))^{p-1}(h_{i+1}(v)(e_i(v)-e_i(u))\\
&=2(u-v)^{p-2}(\ad e_i(u))^{p-2}(h_{i+1}(v)(e_i(v)-e_i(u))^2\\
&=\cdots=p!(h_{i+1}(v)(e_i(v)-e_i(u))^p=0,
\end{align*}
it proves that the case $j=i+1$.
For the case $i=n-1$, if $j=i$ or $j=i+1$, then we can use \eqref{new2}-\eqref{new3} and the same argument to get $(\ad e_i(u))^p(h_j(v))=0$.
When $j=i+2=n+1$, we have $(\ad e_{n-1}(u))^2(h_{n+1}(v))=0$ by \eqref{new22} with $m=1$.
We turn to the case $i=n$ in \eqref{center-eh}.
By \eqref{hiej}, we only have to consider the cases $j\in\{n-1,n,n+1\}$.
If $j=n-1$, then we have $(\ad e_{n}(u))^2(h_{n-1}(v))=0$ by \eqref{new222} with $m=1$.
Moreover, the relation \eqref{new2} implies that $(\ad e_{n}(u))^3(h_{n}(v))=0$.
This proves the case $i=j=n$. When $j=n+1$, a consecutive application of \eqref{new3} and \eqref{en Hmvu}
implies
\begin{align*}
(u-v)^p(\ad e_n(u))^p(h_{n+1}(v))&=p!h_{n+1}(v)(e_n(v)-e_n(u))^p=0,~~for~N=2n,\\
(\ad e_n(u))^p(h_{n+1}(v))&=\mathcal{H}_{v,u}^p=0,~~for~N=2n+1.
\end{align*}

Finally, for \eqref{center-ef}, it follows when $i\neq j$ immediately from \eqref{eifj}.
When $i=j$, we observe using \eqref{new5} repeatedly that
\begin{align*}
(u-v)^{p-1}(\ad e_i(u))^{p-1}(h_{i+1}(v)\tilde{h}_i(v))
=&2\ (u-v)^{p-2}(\ad e_i(u))^{p-2}(h_{i+1}(v)(e_i(v)-e_i(u))\tilde{h}_i(v))\\
=&3! (u-v)^{p-3}(\ad e_i(u))^{p-3}(h_{i+1}(v)(e_i(v)-e_i(u))^2\tilde{h}_i(v))\\
=&\cdots=p!\,h_{i+1}(v)(e_i(v)-e_i(u))^{p-1}\tilde{h}_i(v)=0.
\end{align*}
Hence, $(\ad e_i(u))^{p-1}(h_{i+1}(v)\tilde{h}_i(v))=0$.
We can also set $v=u$ in this identity to see that $(\ad e_i(u))^{p-1}(h_{i+1}(u)\tilde{h}_i(u))=0$.
Then using \eqref{eifj} we obtain for $i\neq n$
\begin{align*}
(u-v)(\ad e_i(u))^p(f_i(v))&=(\ad e_i(u))^{p-1}(k_i(u)-k_i(v))\\
&=(\ad e_i(u))^{p-1}(h_{i+1}(u)\tilde{h}_i(u)-h_{i+1}(v)\tilde{h}_i(v))=0.
\end{align*}
When $i=j=n$ for $N=2n$, we note that $k_n(u)=h_{n+1}(u)\tilde{h}_{n-1}(u)$.
One argues similarly to obtain $(\ad e_n(u))^p(f_n(v))=0$ using \eqref{new7}.
When $i=j=n$ for $N=2n+1$, we have by \eqref{new6} that
\[
(\ad e_n(u))^{p-1}(h_{n+1}(v)\tilde{h}_n(v))=\tilde{\mathcal{H}}_{v,u}^{p-1}=0.
\]
Also, set $v=u$, it gives $ (\ad e_n(u))^{p-1}(h_{n+1}(u)\tilde{h}_n(u))=0$.
Then using again \eqref{eifj} we have for $N=2n+1$ that
\begin{align*}
(u-v)(\ad e_n(u))^p(f_n(v))&=(\ad e_n(u))^{p-1}(k_n(u)-k_n(v))\\
&=(\ad e_n(u))^{p-1}(h_{n+1}(u)\tilde{h}_n(u)-h_{n+1}(v)\tilde{h}_n(v))=0.
\end{align*}
\end{proof}

Now we introduce some useful power series which will be frequently used later on.
Most results follow from \cite[Section 2]{BT18}.

First we assume that $A_I=\cup_{r\geq 0}{\rm F}_r A_I$ a filtered $\kk$-algebra with $1\in{\rm F}_0A_I$ and that we are given elements $\{X^{(r)}\in {\rm F}_{r-1}A_I;~r\geq 0\}$ such that $[\gr_r X^{(r+1)},\gr_s X^{(s+1)}]=0$ for all $r,s\geq 0$.
We necessarily have that $X^{(0)}=0$ since ${\rm F}_{-1}A_I=(0)$ by our conventions.
We define {\it power series of type I}
\[
X_I(u)=\sum\limits_{r\geq 0}X_I^{(r)}u^{-r}:=\big(\sum\limits_{r\geq 0}X^{(r)}u^{-r}\big)^p\in A_I[[u^{-1}]].
\]
Moreover, suppose that $A_{II}$ is a filtered $\kk$-algebra containing commuting elements $\{Y^{(r)}\in {\rm F}_{r-1}A_{II};~r>0\}$.
Setting $Y^{(0)}:=1\in {\rm F}_0 A_{II}$, we define {\it power series of type II}
\[
Y_{II}(u)=\sum\limits_{r\geq 0}Y_{II}^{(r)}u^{-r}:=\prod\limits_{i=1}^p\big(\sum\limits_{r\geq 0}Y^{(r)}(u-i+1)^{-r}\big)\in A_{II}[[u^{-1}]].
\]
\begin{Lemma}\label{useful lemma:power series}\cite[Lemmas 2.1, 2.9]{BT18}
Keep the same notations as above. The following statements hold:
\begin{enumerate}
    \item [(1)] We have that $X_I^{(r)}=0$ for $r<p$ and $X_I^{(p)}=(X^{(1)})^p\in {\rm F}_0 A_I$.
    If $r>p$ and $p\mid r$ then $X_I^{(r)}\in {\rm F}_{r-p}A_I$ and $X_I^{(r)}\equiv (X^{(r/p)})^p ~({\rm mod}~{\rm F}_{r-p-1} A_I)$.
    If $r>p$ and $p\nmid r$ then $X_I^{(r)}\in {\rm F}_{r-p-1}A_I$;
    \item[(2)] We have $Y_{II}^{(0)}=1, Y_{II}^{(r)}=0$ for $1\leq r\leq p-1$, and
    \[
    Y_{II}^{(p)}=(Y^{(1)})^p-Y^{(1)}\in {\rm F}_0 A_{II}.
    \]
   If $r>p$ and $p\mid r$ then $Y_{II}^{(r)}\in {\rm F}_{r-p}A_{II}$ and
   \[
    Y_{II}^{(r)}\equiv (Y^{(r/p)})^p-Y^{(r-p+1)} ~({\rm mod}~{\rm F}_{r-p-1} A_{II}).
   \]
   Finally, if $r>p$ and $p\nmid r$ then $Y_{II}^{(r)}\in {\rm F}_{r-p-1}A_{II}$.
\end{enumerate}
\end{Lemma}

For use in the next theorem, we let
\begin{align}\label{pij-qji}
p_{i,j}(u)=\sum\limits_{r\geq p}p_{i,j}^{(r)}u^{-r}:=e_{i,j}(u)^p,\ \ q_{j,i}(u)=\sum\limits_{r\geq p}q_{j,i}^{(r)}u^{-r}:=f_{j,i}(u)^p.
\end{align}
\begin{Theorem}\label{theorem: off-diag p-center}
For $(i,j)\in I$, the algebras $Z(\xon)\cap\xon_{i,j}^+$ and $Z(\xon)\cap\xon_{j,i}^-$
are infinite rank polynomial algebra freely generated by the central elements $\{p_{i,j}^{(rp)};~r>0\}$ and $\{q_{j,i}^{(rp)};~r>0\}$,
respectively.
We have that $\deg p_{i,j}^{(rp)}=\deg q_{j,i}^{(rp)}=rp-p$ and
\begin{align}\label{greijrp-1}
\gr_{rp-p}p_{i,j}^{(rp)}=(F_{i,j}t^{r-1})^p,~\ \ \gr_{rp-p}q_{j,i}^{(rp)}=(F_{j,i}t^{r-1})^p.
\end{align}
For $r>p$ with $p\nmid r$, we have that $\deg p_{i,j}^{(r)}=r-p-1$ and it is
a polynomial in the elements $\{p_{i,j}^{(sp)};~0<s\leq \lfloor r/p\rfloor\}$.
Analogous statements with $p_{i,j}$ replaced by $q_{j,i}$ also hold.
\end{Theorem}
\begin{proof}
Using the anti-automorphism $\tau$, it suffices to prove all of the statements for $\xon_{i,j}^+$.
We write $A$ for the tensor algebra $\U(\fo_N[t])\otimes \kk[\zeta_1,\zeta_2,\dots]$ and $Z(A)$ for its center.
By Theorem \ref{Thm: center of U}, we know that
\[
Z(A)=Z(\fo_N[t])\otimes\kk[\zeta_1,\zeta_2,\dots]=Z_p(\fo_N[t])\otimes\kk[\zeta_1,\zeta_2,\dots].
\]
Denote by $A_{i,j}^+$ the commutative subalgebra of $A$ generated by $\{F_{i,j}t^r;~r\geq 0\}$.
It follows that
\[
Z(A)\cap A_{i,j}^+=\kk[(F_{i,j}t^r)^p;~r>0].
\]
Note that $\gr(Z(\xon)\cap \xon_{i,j}^+)\subseteq Z(A)\cap A_{i,j}^+$.
By Lemma \ref{offdig-power}, we know that $p_{i,j}^{(rp)}$ belongs to $Z(\xon)\cap \xon_{i,j}^+$.
We apply Lemma \ref{useful lemma:power series}(1) by taking $X^{(r)}:=e_{i,j}^{(r)}$ and $X_I(u)=(e_{i,j}(u))^p=p_{i,j}(u)$ to see that
$\deg p_{i,j}^{(rp)}=rp-p$
and
\[
\gr_{rp-p}p_{i,j}^{(rp)}=(\gr_{r-1}{e_{i,j}^{(r)}})^p=(F_{i,j}t^{r-1})^p.
\]
Consequently, $\gr(Z(\xon)\cap \xon_{i,j}^+)=Z(A)\cap A_{i,j}^+$ and the elements $\{p_{i,j}^{(rp)};~r>0\}$ are algebraically independent generators for $Z(\xon)\cap \xon_{i,j}^+$.
Suppose that $r>p$ with $p\nmid r$.
Applying again Lemma \ref{useful lemma:power series}(1) gives $\deg p_{i,j}^{(r)}=r-p-1$.
Since the coefficient $p_{i,j}^{(r)}$ is central by Lemma \ref{offdig-power},
hence, it is a polynomial in the elements $\{p_{i,j}^{(sp)};~0<s\leq \lfloor r/p\rfloor\}$.
\end{proof}
\begin{Corollary}
For $(i,j)\in I$ and $r>0$, we have that $(e_{i,j}^{(r)})^p,(f_{j,i}^{(r)})^p\in Z(\xon)$ and
\begin{align}\label{greijrp}
\gr_{rp-p}(e_{i,j}^{(r)})^p=(F_{i,j}t^{r-1})^p,~\ \ \gr_{rp-p}(f_{j,i}^{(r)})^p=(F_{j,i}t^{r-1})^p.
\end{align}
\end{Corollary}
\begin{proof}
We only prove the case that $N=2n+1$ and the argument for type $D$ is similar.
Using the anti-automorphism $\tau$ and Lemma \ref{permutation auto},
this reduces to showing that the elements $(e_i^{(r)})^p$ and $(e_{n+1,n+2}^{(r)})^p$ belong to $Z(\xon)$ for all $1\leq i\leq n-1$.
Assume first that $1\leq i\leq n-1$.
Then \eqref{eiei} implies that the elements $e_i^{(r)}$ satisfy the relations
\[
[e_i^{(r)},e_i^{(s)}]=\sum\limits_{t=r}^{s-1}e_i^{(t)}e_i^{(r+s-1-t)}
\]
for all $1\leq r <s$.
According to \cite[Remark 2.2]{BT18},
the elements $p_{i,i+1}^{(rp)}$ can be expressed as polynomials in the elements $\{(e_i^{(s)})^p;~0<s\leq r\}$.
Remember that we already proved $p_{i,i+1}^{(rp)}\in Z(\xon)$.
By induction on $r$ one obtains the elements $\{(e_{i}^{(r)})^p;~r>0\}$ are central in $\xon$.
We recall that $e_{n+1,n+2}(u)=-e_n(u-1/2)$, while $(\alpha_n,\alpha_n)=1$,
the application of \eqref{eiei} with $i=n$ yields
\[
[e_{n+1,n+2}(u),e_{n+1,n+2}(v)]=\frac{1}{2}\frac{\big(e_{n+1,n+2}(v)-e_{n+1,n+2}(u)\big)^2}{u-v}.
\]
We replace $e_{n+1,n+2}(u/2)$ by $\tilde{e}_{n+1,n+2}(u):=\sum_{r\geq 1}\tilde{e}_{n+1,n+2}^{(r)}u^{-r}$.
It follows that
\[
[\tilde{e}_{n+1,n+2}(u),\tilde{e}_{n+1,n+2}(v)]=\frac{\big(\tilde{e}_{n+1,n+2}(v)-\tilde{e}_{n+1,n+2}(u)\big)^2}{u-v}.
\]
By the same token, we have
\[
[\tilde{e}_{n+1,n+2}^{(r)},\tilde{e}_{n+1,n+2}^{(s)}]=\sum\limits_{t=r}^{s-1}\tilde{e}_{n+1,n+2}^{(t)}\tilde{e}_{n+1,n+2}^{(r+s-1-t)}
\]
for all $1\leq r <s$.
The proof of Lemma \ref{offdig-power} also implies that all coefficients in the power series $(e_{n+1,n+2}(u))^p$ belong to $Z(\xon)$.
By the same token, the elements $\{(\tilde{e}_{n+1,n+2}^{(r)})^p;~r>0\}$ are central.
As a result,  $(e_{n+1,n+2}^{(r)})^p\in Z(\xon)$.
\end{proof}
\begin{Remark}
By passing to the associated graded algebra,
we see that the central elements $\{(e_{i,j}^{(r)})^p;~r>0\}$ give another algebraically independent set of generators for $Z(\xon)\cap\xon_{i,j}^+$
lifting the central elements $\{(F_{i,j}t^{r-1})^p;~r>0\}$ of $\gr \xon$.
\end{Remark}

\subsection{Diagonal $p$-central elements}\label{subs: diag-central}
This subsection we introduce the $p$-central elements that belong to the diagonal subalgebras
\[
\xon_i^0:=\kk[h_i^{(r)};~r>0]
\]
of $\xon$.
For $i=1,\dots,n+1$, we define
\begin{align}\label{biu}
 b_i(u)=\sum\limits_{r\geq 0}b_i^{(r)}u^{-r}:=h_{i\downarrow p}(u)=h_i(u)h_i(u-1)\cdots h_i(u-p+1).
\end{align}

\begin{Lemma}\label{lemma:bir in center}
For all $i=1,\dots, n+1$ and $r>0$, the elements $b_i^{(r)}$ belongs to $Z(\xon)$.
\end{Lemma}
\begin{proof}
In view of \eqref{hihj}, it suffices to check that
\[
[b_i(u),e_j(v)=0=[b_i(u),f_j(v)]
\]
for all $1\leq j\leq n$. By applying the anti-automorphism $\tau$ (Subsection \ref{subsec: anto}(1)),
it suffices to check just the first equality.
Assume first that $1\leq i\leq n$.
It is clear when $j\notin\{i-1,i\}$ by \eqref{hiej}.
 When $j=i$ or $j=i-1$, the identities \eqref{down 1}-\eqref{up 1} imply that
\[
  [b_i(u),e_{i}(v)]=[h_{i\downarrow p}(u),e_{i}(v)]=0
 \]
  and
 \[
 [b_i(u),e_{i-1}(v)]=[h_{i\downarrow p}(u),e_{i-1}(v)]=[h_{i\uparrow p}(u-p+1),e_{i-1}(v)]=0.
 \]
If $i=n+1$ and $N=2n$,  then \eqref{hiej} again implies that we may assume that $j=n$ or $j=n-1$.
One argues similarly for this case using \eqref{up 2-2n}-\eqref{down 2-2n} instead of  \eqref{down 1}-\eqref{up 1}.
We now treat the case $i=n+1$ and $N=2n+1$.
In view of Theorem \ref{thm: HC center via Drinfeld} and \eqref{hihj},
we obtain
\begin{align*}
c_{\downarrow p}(u)&:=c(u)c(u-1)\cdots c(u-p+1)\\
&=h_{n+1\downarrow p}(u-n+1/2)h_{n+1\downarrow p}(u-n)\\
&=b_{n+1}(u-n+1/2)b_{n+1}(u-n)=(b_{n+1}(u))^2.
\end{align*}
Since all the coefficients of $c(u)$ are central, so are $c_{\downarrow p}(u)$,
the centrality of $b_{n+1}^{(r)}$ follows by induction.
\end{proof}

\begin{Theorem}\label{thm: diag center}
For $1\leq i\leq n$,
the algebra $Z(\xon)\cap\xon_i^0$ is an infinite rank polynomial algebra freely generated by the central elements $\{b_i^{(rp)};~r>0\}$.
The statement also holds for $i=n+1$ and $N=2n$.
We have that $\deg b_i^{(rp)}=rp-p$ and
\begin{align}\label{gr birp}
\gr_{rp-p}b_i^{(rp)}=(F_{i,i}t^{r-1})^p-F_{i,i}t^{rp-p}+(\zeta_{r}^p-\zeta_{rp-p+1})/2.
\end{align}
For $1<r<p$, we have that $b_i^{(r)}=0$.
For $r>p$ with $p\nmid r$, we have that $\deg b_i^{(r)}=r-p-1$ and it is
a polynomial in the elements $\{b_i^{(sp)};~0<s\leq \lfloor r/p\rfloor\}$.
\end{Theorem}
\begin{proof}
Similar to the proof of Theorem \ref{theorem: off-diag p-center},
we denote by $A$ the tensor algebra $\U(\fo_N[t])\otimes \kk[\zeta_1,\zeta_2,\dots]$ and $Z(A)$ for its center.
Let $A_i^0$ be the commutative subalgebra of $A$ generated by $\{F_{i,i}t^r+\zeta_{r+1}/2;~r\geq 0\}$.
One sees that
\[
Z(A)\cap A_i^0=\kk[(F_{i,i}t^r+\zeta_{r+1}/2)^p-(F_{i,i}t^{rp}+\zeta_{rp+1}/2);~r\geq 0].
\]
We have that $\gr(Z(\xon)\cap\xon_i^0)\subset Z(A)\cap A_i^0$.
By Lemma \ref{lemma:bir in center}, we know that $b_i^{(rp+p)}$ belongs to $Z(\xon)\cap\xon_i^0$.
Moreover, applying Lemma \ref{useful lemma:power series}(2) with $Y^{(r)}:=h_i^{(r)}$ and $Y_{II}(u)=b_i(u)$,
we see that $\deg b_i^{(rp+p)}=rp$ and
\[
\gr_{rp}b_{i}^{(rp+p)}=(F_{i,i}t^{r}+\zeta_{r+1}/2)^p-(F_{i,i}t^{rp}+\zeta_{rp+1}/2).
\]
We thus obtain $\gr(Z(\xon)\cap\xon_i^0)=Z(A)\cap A_i^0$ and the elements $\{b_i^{(rp)};~r>0\}$ are algebraically independent generators for $Z(\xon)\cap\xon_i^0$.
Again, Lemma \ref{useful lemma:power series}(2) implies that $b_i^{(r)}=0$ for $1<r<p$ and that $\deg b_i^{(r)}=r-p-1$ if $r>p$ with $p\nmid r$.
Since it is central by Lemma \ref{lemma:bir in center},
it must be a polynomial in the elements $\{b_i^{(sp)};~0<s\leq \lfloor r/p\rfloor\}$.
\end{proof}
\subsection{The center $Z(\xon)$}
We define the {\it $p$-center} $Z_p(\xon)$ of $\xon$ to be the subalgebra generated by
\begin{align}\label{p-center generators}
\{b_i^{(rp)};~1\leq i\leq n+1, r>0\} \cup \big\{\big(e_{i,j}^{(r)}\big)^p, \big(f_{j,i}^{(r)}\big)^p;~(i,j)\in I, r>0\big\}
\end{align}
According to Lemma \ref{offdig-power} and Lemma \ref{lemma:bir in center},
we know that both $Z_{\HC}(\xon)$ and $Z_p(\xon)$ are subalgebras of $Z(\xon)$.

We let
\begin{align*}
bc(u):=\sum\limits_{r\geq 0} bc^{(r)}u^{-r}:=c(u)c(u-1)\cdots c(u-p+1).
\end{align*}
Then Theorem \ref{thm: HC center via Drinfeld} readily implies
\begin{align}\label{bcu}
bc(u)=\left\{
\begin{array}{ll}
(b_{n+1}(u))^2 &  for~N=2n+1,\\
b_n(u)b_{n+1}(u)  &  for~N=2n.
\end{array}
\right.
\end{align}
For example, let us verify \eqref{bcu} in the case $N=2n+1$.
In view of Theorem \ref{thm: HC center via Drinfeld},
we have
\begin{align*}
bc(u)&=\prod\limits_{i=1}^p c(u-i+1)=\prod\limits_{i=1}^p \prod\limits_{j=1}^n\frac{h_j(u-i-j+2)}{h_j(u-i-j+1)}\times h_{n+1}(u-i-n+3/2)h_{n+1}(u-i-n+1)\\
&=\prod\limits_{i=1}^p\frac{h_1(u-i+1)}{h_1(u-i)}\prod\limits_{i=1}^p\frac{h_2(u-i)}{h_2(u-i-1)}\times\cdots \times\prod\limits_{i=1}^p\frac{h_n(u-i-n+2)}{h_n(u-i-n+1)}\\
 &\quad\quad \times \prod\limits_{i=1}^p h_{n+1}(u-i-n+3/2)\times \prod\limits_{i=1}^p h_{n+1}(u-i-n+1)\\
 &=1\cdot b_{n+1}(u)\cdot b_{n+1}(u)=(b_{n+1}(u))^2.
\end{align*}

By definition each $bc^{(r)}$ can be expressed as a polynomial in the elements $\{c^{(s)};~s>0\}$,
so that it belongs to $Z_{\HC}(\xon)$.
It is also a polynomial in the elements $\{b_i^{(r)};~i=n+1, (n+1)', r>0\}$,
so that it belongs to $Z_p(\xon)$ by Theorem \ref{thm: diag center}.
Consequently, $bc^{(r)}\in Z_{\HC}(\xon)\cap Z_p(\xon)$.
\begin{Lemma}
For $r>0$, we have $\deg bc^{(rp)}=rp-p$ and
\begin{align}\label{gr bcrp}
\gr_{rp-p} bc^{(rp)}=\zeta_r^p-\zeta_{rp-p+1}.
\end{align}
\end{Lemma}
\begin{proof}
We just go through the case $N=2n+1$, since another case is similar.
Set
\[H_{n+1}(u):=1+\sum_{r>0}H_{n+1}^{(r)}u^{-r}:=(h_{n+1}(u))^2.\]
It is easy to check that $\deg H_{n+1}^{(r)}=r-1$ and $\gr_{r-1}H_{n+1}^{(r)}=2\gr_{r-1}h_{n+1}^{(r)}$.
Observe that
\[bc(u)=H_{n+1}(u)H_{n+1}(u-1)\cdots H_{n+1}(u-p+1).\]
Applying Lemma \ref{useful lemma:power series}(2) with $Y^{(r)}:=H_{n+1}^{(r)}$ and $Y_{II}(u)=bc(u)$,
we have
$\deg bc^{(rp)}=rp-p$ and
\[
\gr_{rp-p}bc^{(rp)}=(\gr_{r-1}H_{n+1}^{(r)})^p-\gr_{rp-p}H_{n+1}^{(rp-p+1)}=(2\gr_{r-1}h_{n+1}^{(r)})^p-2\gr_{rp-p}h_{n+1}^{(rp-p+1)}.
\]
Our assertion now follows from \eqref{identification-2}.
\end{proof}

\begin{Theorem}\label{main theorem1: centre of  Xon}
The center $Z(\xon)$ is generated by $Z_{\HC}(\xon)$ and
$Z_p(\xon)$. Moreover:
\begin{enumerate}
\item $Z_{\HC}(\xon)$ is the free polynomial algebra generated by $\{c^{(r)};~r > 0\}$;
\item $Z_p(\xon)$ is the free polynomial algebra generated by
\begin{equation}\label{generator of p centre}
\{b_i^{(rp)};~1\leq i\leq n+1, r>0\} \cup \big\{\big(e_{i,j}^{(r)}\big)^p, \big(f_{j,i}^{(r)}\big)^p;~(i,j)\in I, r>0\big\};
\end{equation}
\item $Z(\xon)$ is the free polynomial algebra generated by
\begin{equation}\label{generator of centre }
\{b_i^{(rp)}, c^{(r)};~1\leq i\leq n, r>0\} \cup \big\{\big(e_{i,j}^{(r)}\big)^p, \big(f_{j,i}^{(r)}\big)^p;~(i,j)\in I, r>0\big\};
\end{equation}
\item $Z_{\HC}(\xon)\cap Z_p(\xon)$ is the free polynomial algebra generated by $\{bc^{(rp)};~r>0\}$.
\end{enumerate}
\end{Theorem}
\begin{proof}
(1) This was proved in \cite[Corollary 3.9]{AMR06}, see also \cite{AACFR03}.

(2) The given elements generate $Z_p(\xon)$ by the definition.
It follows from \eqref{greijrp} and \eqref{gr birp} that they are lifts of the algebraically independent elements of the
associated graded algebra $\U(\fo_N[t])\otimes \kk[\zeta_1,\zeta_2,\dots]$.
This yields the assertion.

(3) Let $Z$ be the subalgebra of $Z(\xon)$ generated by the given elements.
In view of Theorem \ref{Thm: center of U}, we have
\[
\gr Z\subseteq \gr Z(\xon)\subseteq Z(\gr \xon)=Z_p(\fo_N[t])\otimes\kk[\zeta_1,\zeta_2,\dots].
\]
Then the foregoing observations in conjunction with Theorem \ref{Thm: center of U} imply
 that the generators of $Z$ are lifts of the algebraically independent elements of $Z_p(\fo_N[t])\otimes\kk[\zeta_1,\zeta_2,\dots]$.
Moreover,  by \eqref{greijrp} and \eqref{gr birp} , it is clear that $Z_p(\fo_N[t])\otimes\kk[\zeta_1,\zeta_2,\dots]\subseteq \gr Z$.
This implies that $Z=Z(\xon)$.

(4) We already observed that all $bc^{(rp)}$ belong to $Z_{\HC}(\xon)\cap Z_p(\xon)$,
while \eqref{gr bcrp} implies they are algebraically independent.
We claim that $Z_p(\xon)$ is freely generated by the elements
\[
\{b_i^{(rp)}, bc^{(rp)};~1\leq i\leq n, r>0\} \cup \big\{\big(e_{i,j}^{(r)}\big)^p, \big(f_{j,i}^{(r)}\big)^p;~(i,j)\in I, r>0\big\}.
\]
We know already from (3) that all of these elements different from $bc^{(rp)}$ are algebraically independent of anything in $Z_{\HC}(\xon)$.
Our result thus follows from the claim. To prove the claim, we use \eqref{greijrp}, \eqref{gr birp} and \eqref{gr bcrp} to pass to the associated graded algebra.
Let $A^0$ be the subalgebra of $\U(\fo_N[t])\otimes \kk[\zeta_1,\zeta_2,\dots]$ generated by $\{\gr_{rp-p} b_i^{(rp)};~1\leq i\leq n+1, r>0\}$.
By (2), it suffices to verify that
\[
\{\gr_{rp-p} b_i^{(rp)}, \gr_{rp-p} bc^{(rp)};~1\leq i\leq n, r>0\}
\]
freely genetrate $A^0$. This is easily seen from \eqref{gr birp} and \eqref{gr bcrp}.
\end{proof}

\subsection{The center $Z(\yon)$}
The Yangian $\yon$ is defined as the subalgebra of $\xon$ which consists of the elements stable under the automorphisms
\begin{align}\label{auto:mu_f}
\mu_f:T(u)\rightarrow f(u)T(u)
\end{align}
for all series $f(u)\in 1+u^{-1}\kk[[u^{-1}]]$. That is,
\[
\yon:=\{x\in\xon;~\mu_f(x)=x~{\rm for~all~}f(u)\in 1+u^{-1}\kk[[u^{-1}]]\}.
\]
By the same proof as \cite[Theorem 3.1]{AMR06}, we have the tensor product decomposition
\begin{align}\label{tensor product decom xon}
\xon=Z_{\HC}(\xon)\otimes \yon.
\end{align}
On Gaussian generators, clearly
\[
\mu_f(h_i(u))=f(u)h_i(u),~\mu_f(e_{i,j}(u))=e_{i,j}(u)~\text{and}~\mu_f(f_{j,i}(u))=f_{j,i}(u).
\]

Recall from \eqref{def:kef-1} and \eqref{def:kn for D} that
\[
k_i(u)=\sum\limits_{r\geq 0}k_i^{(r)}u^{-r}:=\left\{
\begin{array}{ll}
\tilde{h}_i(u)h_{i+1}(u),&~if~1\leq i\leq n-1,\\
\tilde{h}_n(u)h_{n+1}(u),&~if~i=n~and~N=2n+1,\\
\tilde{h}_{n-1}(u)h_{n+1}(u),&~if~i=n~and~N=2n.
\end{array}
\right.
\]
Since $\tilde{h}_i(u)=\sum_{r\geq 0}\tilde{h}_i^{(r)}u^{-r}=h_i(u)^{-1}$,
we have $\tilde{h}_i^{(0)}=1$ and $\tilde{h}_{i}^{(r)}=-\sum_{t=1}^{r}h_i^{(t)}\tilde{h}_i^{(r-t)}$.
In particular, $\deg k_i^{(r+1)}=r$ and $\gr_r h_i^{(r+1)}=-\gr_r \tilde{h}_i^{(r+1)}$.
Moreover, the identifications \eqref{identification-1}-\eqref{identification-2} yield:
\begin{align}\label{grkir}
\gr_r k_i^{(r+1)}=\left\{
\begin{array}{ll}
F_{i+1,i+1}t^r-F_{i,i}t^r,&~if~1\leq i\leq n-1,\\
-F_{n,n}t^r,&~if~i=n~and~N=2n+1,\\
-F_{n,n}t^r-F_{n-1,n-1}t^r,&~if~i=n~and~N=2n.
\end{array}
\right.
\end{align}
The decomposition \eqref{tensor product decom xon} and Theorem \ref{thm: HC center via Drinfeld} imply that the subalgebra $\yon$ of $\xon$ is generated by
\[
\{k_i^{(r)};~1\leq i\leq n, r>0\}\cup \{e_i^{(r)}, f_i^{(r)};~1\leq i\leq n, r>0\},
\]
see \cite[Proposition 6.1]{JLM18}.
The ascending filtration on the Yangian $\yon$ is induced by the one on $\xon$.
In particular, under the identification \eqref{PBW map}, we have $\gr \yon=\U(\fo_N[t])$ (cf. \cite[Theorem 3.6]{AMR06}).

Let
\begin{align}\label{aiu}
a_i(u)&=\sum_{r \geq 0}a_i^{(r)} u^{-r}:=k_i(u)k_i(u-1)\cdots k_i(u-p+1)\nonumber\\
&=\left\{
\begin{array}{ll}
b_{i+1}(u)b_i(u)^{-1}&~if~1\leq i\leq n-1,\\
b_{n+1}(u)b_n(u)^{-1}&~if ~i=n~and~N=2n+1,\\
b_{n+1}(u)b_{n-1}(u)^{-1}&~if~i=n~and~N=2n,
\end{array}\right.
\end{align}
where the last equality follows from the definitions (\eqref{def:kef-1}, \eqref{def:kn for D} and \eqref{biu}).
In view of Lemma \ref{lemma:bir in center},
each $a_i^{(r)}$ belongs to $Z(\yon)$.
We define the $p$-center of $\yon$ to be the subalgebra $Z_p(\yon)$ of $Z(\yon)$ generated by
\begin{align}\label{generators of zpyon}
\{a_i^{(rp)};~1\leq i\leq n, r>0\} \cup \big\{\big(e_{i,j}^{(r)}\big)^p, \big(f_{j,i}^{(r)}\big)^p;~(i,j)\in I, r>0\big\}.
\end{align}

\begin{Theorem}\label{main theorem2: centre of  Xon}
The generators \eqref{generators of zpyon} of $Z_p(\yon)$ are algebraically independent, and we have that $\gr Z_p(\yon)=Z_p(\fo_N[t])=Z(\fo_N[t])$.
In particular, $Z_p(\yon)=Z(\yon)$.
\end{Theorem}
\begin{proof}
From \eqref{aiu} , we know that $\deg a_i^{(rp)}=rp-p$ and
\[
\gr_{rp-p} a_i^{(rp)}=\left\{
\begin{array}{ll}
\gr_{rp-p}  b_{i+1}^{(rp)}-\gr_{rp-p}b_i^{(rp)}&~if~1\leq i\leq n-1,\\
\gr_{rp-p}  b_{n+1}^{(rp)}-\gr_{rp-p}b_n^{(rp)}&~if~i=n~and~N=2n+1,\\
\gr_{rp-p}  b_{n+1}^{(rp)}-\gr_{rp-p}b_{n-1}^{(rp)}&~if~i=n~and~N=2n.
\end{array}\right.
\]
It follows from \eqref{gr birp} that
\[
\gr_{rp-p} a_i^{(rp)}=\left\{
\begin{array}{ll}
\!\!\!(F_{i+1,i+1}t^{r-1}-F_{i,i}t^{r-1})^p-(F_{i+1,i+1}t^{rp-p}-F_{i,i}t^{rp-p})&\!\!\!\!if~1\leq i\leq n-1,\\
\!\!\!F_{n,n}t^{rp-p}-(F_{n,n}t^{r-1})^p&\!\!\!\!if~i=n~and~N=2n+1,\\
\!\!\!(F_{n,n}t^{rp-p}+F_{n-1,n-1}t^{rp-p})-(F_{n,n}t^{r-1}+F_{n-1,n-1}t^{r-1})^p&\!\!\!\!if~i=n~and~N=2n.
\end{array}\right.
\]
Combining with \eqref{greijrp}, we see that the generators \eqref{generators of zpyon} are lifts of generator for $Z_p(\fo_N[t])$.
This establishes the algebraic independence and that $\gr Z_p(\yon)=Z_p(\fo_N[t])$,
and the remaining assertion thus follows from Theorem \ref{Thm: center of U}.
\end{proof}

\bigskip
\appendix
\section{Proof of the Serre relations}\label{section name:appendix}
In the appendix,
we will prove the Serre relations (\eqref{Serre-e} and \eqref{Serre-f}).
We establish this result by adapting the methods of \cite{Lev93, GNW18} to positive characteristic.

Define the series by
\begin{align*}
\kappa_i(u):=h_i(u-(i-1)/2)^{-1}h_{i+1}(u-(i-1)/2)
\end{align*}
for $i=1,\dots,n-1$, and
\[
\kappa_n(u):=\left\{
\begin{array}{ll}
h_n(u-(n-1)/2)^{-1}h_{n+1}(u-(n-1)/2)   &~\text{for}~N=2n+1,\\
h_{n-1}(u-(n-2)/2)^{-1}h_{n+1}(u-(n-2)/2)   & ~\text{for}~N=2n.
\end{array}{}
\right.
\]
Moreover, set
\[
    \xi_i^+(u):=f_{i+1,i}(u-(i-1)/2),~\quad \xi_i^-(u):=e_{i,i+1}(u-(i-1)/2)
\]
for $i=1,\dots,n-1$,
\[
\xi_n^+(u):=\left\{
\begin{array}{ll}
f_{n+1,n}(u-(n-1)/2) &~\text{for}~N=2n+1,\\
f_{n+1,n-1}(u-(n-2)/2)  & ~\text{for}~N=2n,
\end{array}{}
\right.
\]
and
\[
\xi_n^-(u):=\left\{
\begin{array}{ll}
e_{n,n+1}(u-(n-1)/2) &~\text{for}~N=2n+1,\\
e_{n-1,n+1}(u-(n-2)/2)  & ~\text{for}~N=2n,
\end{array}{}
\right.
\]
Introduce elements by the respective expansions into power series in $u^{-1}$,
\[
\kappa_i(u)=1+\sum\limits_{r\geq 0}\kappa_{i,r}u^{-r-1}
\quad \text{and}\quad \xi_i^{\pm}(u)=\sum\limits_{r\geq 0}\xi_{i,r}^{\pm}u^{-r-1}
\]
for $i=1,\dots,n$.

\begin{Lemma}\label{Lemma:app-1}
The following relations holds in $\xon$:
\begin{align}\label{D yangian rel-1}
[\kappa_{i,r},\kappa_{j,s}]=0,
\end{align}
\begin{align}\label{D yangian rel-2}
[\xi_{i,r}^+,\xi_{j,s}^-]=\delta_{i,j}\kappa_{i,r+s},
\end{align}
\begin{align}\label{D yangian rel-3}
[\kappa_{i,0},\xi_{j,s}^{\pm}]=\pm(\alpha_i,\alpha_j)\xi_{j,s}^{\pm},
\end{align}
\begin{align}\label{D yangian rel-4}
[\kappa_{i,r+1},\xi_{j,s}^{\pm}]-[\kappa_{i,r},\xi_{j,s+1}^{\pm}]=\pm\frac{(\alpha_i,\alpha_j)}{2}(\kappa_{i,r}\xi_{j,s}^{\pm}+\xi_{j,s}^{\pm}\kappa_{i,r}),
\end{align}
\begin{align}\label{D yangian rel-5}
[\xi_{i,r+1}^{\pm},\xi_{j,s}^{\pm}]-[\xi_{i,r}^{\pm},\xi_{j,s+1}^{\pm}]=\pm\frac{(\alpha_i,\alpha_j)}{2}(\xi_{i,r}\xi_{j,s}^{\pm}+\xi_{j,s}^{\pm}\xi_{i,r}).
\end{align}
\end{Lemma}
\begin{proof}
The proof amounts to rephrasing the relations \eqref{hihj}-\eqref{fifj} of Theorem \ref{Theorem: Drinfeld presentation} in terms of the series $\kappa_i(u)$ and $\xi_{i}^{\pm}(u)$, see the proof \cite[Proposition 6.2]{JLM18}.
\end{proof}

Inspired by the work of \cite{Lev93},
we introduce the elements $\tka_{i,r}$ for $r=1,2,3$ via
\begin{align*}
\tka_{i,1}:=\kappa_{i,1}-\frac{1}{2}\kappa_{i,0}^2, \quad \tka_{i,2}:=\kappa_{i,2}-\kappa_{i,0}\kappa_{i,1}+\frac{1}{3}\kappa_{i,0}^3,
\end{align*}
\begin{align*}
\tka_{i,3}:=\kappa_{i,3}-\kappa_{i,0}\kappa_{i,2}-\frac{1}{2}\kappa_{i,1}^2+\kappa_{i,0}^2\kappa_{i,1}-\frac{1}{4}\kappa_{i,0}^4.
\end{align*}

The following Lemma is a special case of \cite[Lemma 1.4]{Lev93},
the proof given in characteristic zero in {\it loc.cit.} (see also \cite[Section 2.9]{GTL13}).

\begin{Lemma}\label{Lemma:app-2}
The following relations are satisfied for all $1\leq i,j\leq n$ and $s\geq 0$:
 \begin{align}\label{tka xjs bracket-1}
 [\tka_{i,1},\xi_{j,s}^{\pm}]=\pm(\alpha_i,\alpha_j)\xi_{j,s+1}^{\pm},\quad
[\tka_{i,2},\xi_{j,s}^{\pm}]=\pm(\alpha_i,\alpha_j)\xi_{j,s+2}^{\pm}\pm\frac{1}{12}(\alpha_i,\alpha_j)^3\xi_{j,s}^{\pm},
 \end{align}
\begin{align}\label{tka xjs bracket-2}
[\tka_{i,3},\xi_{j,s}^{\pm}]=\pm(\alpha_i,\alpha_j)\xi_{j,s+3}^{\pm}\pm \frac{1}{4}(\alpha_i,\alpha_j)^3\xi_{j,s+1}^{\pm}.
 \end{align}
\end{Lemma}
\begin{proof}
These follow from \eqref{D yangian rel-3} and \eqref{D yangian rel-4} by a direct computation.
For example, the case $s=0$ in the second identity \eqref{tka xjs bracket-1}
has been proved in \cite[Proposition 2.31]{GNW18}.
We can repeat the argument there without any changes to the general case.
\end{proof}

Now we are ready to prove the Serre relations.
As noticed in the proof of \cite[Theorem 5.14]{JLM18},
it suffices to prove the following Proposition,
which implies the Serre relations in the algebra $\xon$.
\begin{Proposition}\label{prop:app-serre relation}
For any $i\neq j$, the following identities hold:
\begin{align}\label{D serre rel}
\sum_{\sigma\in\fS_k}\big[\xi_{i,r_{\sigma(1)}}^{\pm},
\big[\xi_{i,r_{\sigma(2)}}^{\pm},\dots,\big[\xi_{i,r_{\sigma(k)}}^{\pm},\xi_{j,s}^{\pm}\big]\dots\big]\big]=0,
\end{align}
for $k=1+|c_{i,j}|$.
\end{Proposition}
\begin{proof}
We only prove \eqref{D serre rel} for the case $i=n$ and $j=n-1$ with $N=2n+1$,
as others can be treated similarly.
Note that $c_{n,n-1}=-2$.
In this case, denote the left-hand side of \eqref{D serre rel} by $\Qpm(r_1,r_2,r_3;s)$.
Our goal is to prove
$\Qpm(r_1,r_2,r_3;s)=0$ for any $r_1,r_2,r_3,s\geq 0$.
Note that $\Qpm$ symmetric with respect to the variables $r_1,r_2,r_3$.
We first show that $\Qpm(0,0,0;s)=0$ by induction on $s\geq 0$.
This has already been proven in \cite[Lemma 2.33]{GNW18}.
If $s=0$,
then this is implied by the embedding $\U(\fo_N)\hookrightarrow \xon$ in combination with the Serre relations in the Lie algebra $\fo_N$.
Suppose that $\Qpm(0,0,0;s)=0$ for some $s\geq 0$.
Thanks to Lemma \ref{Lemma:app-2}, we may apply $[\tka_{n,1},\cdot]$ and $[\tka_{n-1,1},\cdot]$ to $\Qpm(0,0,0;s)$ to get
\begin{align*}
3\Qpm(1,0,0;s)&-\Qpm(0,0,0;s+1)=0\\
-3\Qpm(1,0,0;s)&+2\Qpm(0,0,0;s+1)=0.
\end{align*}
It follows that
\[
\Qpm(1,0,0;s)=0=\Qpm(0,0,0;s+1).
\]
Therefore, by induction we have $\Qpm(0,0,0;s)=0$ for all $s\geq 0$.
We simultaneously have proven that $\Qpm(1,0,0;s)=0$ for all $s\geq 0$.
Next, we prove the following claim.

{\it Claim.} Suppose that $\Qpm(r_1,r_2,r_3;s)=0$ for all $s\geq 0$.
Then we have
\[
\Qpm(r_1+i,r_2,r_3;s)+\Qpm(r_1,r_2+i,r_3;s)+\Qpm(r_1,r_2,r_3+i;s)=0
\]
for all $s\geq 0$ and $i=1,2,3$.

We prove the claim by applying again Lemma \ref{Lemma:app-2}.
We apply $[\tka_{n,1},\cdot]$ to $\Qpm(r_1,r_2,r_3;s)=0$.
By \eqref{tka xjs bracket-1} we have
\begin{align}\label{app-proof-1}
\Qpm(r_1+1,r_2,r_3;s)+\Qpm(r_1,r_2+1,r_3;s)+\Qpm(r_1,r_2,r_3+1;s)=0.
\end{align}
This proves the claim for $i=1$.
Next, we apply $[\tka_{n,2},\cdot]$ to $\Qpm(r_1,r_2,r_3;s)=0$. We have
\begin{align}\label{app-proof-2}
0=\Qpm(r_1+2&,r_2,r_3;s)+\Qpm(r_1,r_2+2,r_3;s)+\Qpm(r_1,r_2,r_3+2;s)\nonumber \\
-&\frac{1}{6}\Qpm(r_1,r_2,r_3;s)-\Qpm(r_1,r_2,r_3;s+2).
\end{align}
Since the last two terms vanish, the claim for $i=2$ follows.
We apply $[\tka_{n,3},\cdot]$ to $\Qpm(r_1,r_2,r_3;s)=0$.
By \eqref{tka xjs bracket-2} we have
\begin{align}\label{app-proof-3}
0=\Qpm(r_1+3,r_2,&r_3;s)+\Qpm(r_1,r_2+3,r_3;s)+\Qpm(r_1,r_2,r_3+3;s)\nonumber \\
+\frac{1}{4}(\Qpm&(r_1+1,r_2,r_3;s)+\Qpm(r_1,r_2+1,r_3;s)+\Qpm(r_1,r_2,r_3+1;s))\nonumber \\
-&\Qpm(r_1,r_2,r_3;s+3)-\frac{1}{4}\Qpm(r_1,r_2,r_3;s+1).
\end{align}
The assumption in conjunction with \eqref{app-proof-1} implies that the claim is true for $i=3$.

Since we have proven that $\Qpm(0,0,0;s)=\Qpm(1,0,0;s)=0$ for all $s\geq 0$,
it follows immediately that
\[
\Qpm(2,0,0;s)=\Qpm(3,0,0;s)=0
\]
and
\begin{align*}
\Qpm(2,0,0;s)+\Qpm(1,1,0;s)+\Qpm(1,0,1;s)=\Qpm(2,0,0;s)+2\Qpm(1,1,0;s)=0\\
\Qpm(3,0,0;s)+\Qpm(1,2,0;s)+\Qpm(1,0,2;s)=\Qpm(3,0,0;s)+2\Qpm(2,1,0;s)=0
\end{align*}
As a result, we obtain $\Qpm(1,1,0;s)=\Qpm(2,1,0;s)=0$ for all $s\geq 0$.
Moreover, this yields
\[
 2\Qpm(2,1,0;s)+\Qpm(1,1,1;s)=0,
\]
whence $\Qpm(1,1,1;s)=0$ for all $s\geq 0$.
So far we have proved $\Qpm(r_1,r_2,r_3;s)=0$ for all $r_1+r_2+r_3\leq 3$ and $s\geq 0$.
By induction on $r_1+r_2+r_3$, one can derive the desired equality $\Qpm(r_1,r_2,r_3;s)=0$ for any $r_1,r_2,r_3,s\geq 0$.
\end{proof}

\noindent
\textbf{Acknowledgment.}
We would like to thank Ming Liu, Kang Lu and Alexander Molev for helpful discussions.
We are particularly grateful to the referee for his or her helpful comments and suggestions. This work is supported by the National Natural
Science Foundation of China (Grant Nos. 11801394, 12461005),
the Fundamental Research Funds for the Central Universities (No. CCNU25JCPT031) and
the Natural Science Foundation of Hubei Province (No. 2025AFB716).

\end{document}